\documentclass[10pt]{amsart}
\usepackage{amsmath}
\usepackage{amsthm}
\usepackage{amsfonts}
\usepackage{amssymb}
\usepackage{latexsym}
\usepackage{enumerate}
\usepackage{color}

\newcommand{\beq}{\begin{eqnarray}}
\newcommand{\eeq}{\end{eqnarray}}
\newcommand{\beqno}{\begin{eqnarray*}}
\newcommand{\eeqno}{\end{eqnarray*}}
\newcommand{\bR}{\mathbb R}
\newcommand{\bN}{\mathbb N}
\newcommand{\cN}{\mathcal N}

\newcommand{\Ruzicka}{R\r{u}\v zi\v cka}

\newcommand{\midint}{- \mskip-18mu \int}
\newcommand{\Hom}{\mbox{Hom}}

\renewcommand{\phi}{\varphi}

\newcommand {\s}{\tilde{s}}

\DeclareMathOperator{\spt}{supp}

\allowdisplaybreaks

\def\mvint_#1{\mathchoice
  {\mathop{\vrule width 6pt height 3 pt depth -2.5pt
        \kern -8pt \intop}\nolimits_{\kern -3pt #1}}%
 {\mathop{\vrule width 5pt height 3 pt depth -2.6pt
                  \kern -6pt \intop}\nolimits_{#1}}%
 {\mathop{\vrule width 5pt height 3 pt depth -2.6pt
                  \kern -6pt \intop}\nolimits_{#1}}%
 {\mathop{\vrule width 5pt height 3 pt depth -2.6pt
                  \kern -6pt \intop}\nolimits_{#1}}}

\newtheorem{theorem}{Theorem}[section]

\newtheorem{Lem}[theorem]{Lemma}
\newtheorem{corollary}[theorem]{Corollary}

\newtheorem*{Bem}{Remark}

\theoremstyle{definition}
\newtheorem{definition}[theorem]{Definition}
\newtheorem*{remark}{Remark}

\numberwithin{equation}{section}

\begin{document}

\title{Calder\'on--Zygmund estimates for higher order systems with $p(x)$ growth}

\author[J. Habermann]{Jens Habermann}
\address{Jens Habermann, Institute for mathematics\\
Friedrich-Alexander University\\ Bismarckstr. 1 1/2\\ 91054
Erlangen\\ Germany;}
\email{habermann@mi.uni-erlangen.de}

\begin{abstract}
For weak solutions $u \in W^{m,1}(\Omega;\bR^N)$ of higher order systems of the type 
\begin{equation*}
     \int_\Omega \left< A(x,D^m u),D^m \phi \right>\, dx = \int_\Omega \left< |F|^{
     p(x)-2}F,D^m \phi \right>\, dx, \quad \mbox{for all } \phi \in C^{\infty}_c(\Omega;\bR^N), \quad m > 1,
\end{equation*}
with variable growth exponent $p:\Omega \to (1,\infty)$ we prove that if $|F|^{p(\cdot)} \in L^{q}_{loc}(\Omega)$ with $ 1 < q < 
\frac{n}{n-2} + \delta$, then $|D^m u|^{p(\cdot)} \in L^q_{loc}(\Omega)$. We should note that we prove this implication both in the non--degenerate ($\mu > 0$) and in the degenerate case ($\mu = 0$).
\end{abstract}

\maketitle

\section{Introduction}

In this paper we are concerned with a regularity result for weak solutions of
systems of higher order with $p(x)$-- growth. 

Let $n \in \bN_{\ge 2}, N \in \bN_{\ge 1}$ 
and $\Omega \subset \bR^n$ a bounded domain. We consider weak solutions of the 
system
\begin{equation}\label{Gleichung}
     \int_{\Omega} \left<A\left(x,D^m u\right),D^m \phi\right> \, dx = \int_{\Omega} 
     \bigl<|F(x)|^{p(x)-2}F(x),D^m \phi\bigr>\, dx,
\end{equation}
for $\phi \in W^{m,1}_0\left(\Omega;\bR^N\right)$ with $\left|D^m \phi\right|^{p
(\cdot)} \in L^1_{loc} \left(\Omega\right),\ \spt \phi \Subset \Omega$.
Here $A$ denotes a vector field $A: \Omega \times \odot^m(\bR^n,\bR^N) \to 
\Hom(\odot^m(\bR^n,\bR^N),\bR)$, $F: \Omega \to 
\bR^{N\binom{n+m-1}{m}}$, and $p: \Omega \to (1,\infty)$ a measurable function. 
$\odot^m(\bR^n,\bR^N)$ denotes the space of symmetric $m$-- linear forms on $\bR^n$
with values in $\bR^N$. 
The coefficient $A$ is supposed to have $p(x)$-- growth, i.e. for $\mu \in [0,1]$ 
there holds
\begin{equation*}
     \left< D_z A\left(x,z\right) \lambda, \lambda\right> \approx
     \left(\mu^2+ |z|^2\right)^{\frac{p(x)-2}{2}}|\lambda|^2.
\end{equation*}
Additionally we assume that the coefficient $A$ is continuous with respect to the 
first variable and that there exists a modolus of continuity for the exponent
function $p$, which satisfies
\begin{equation}\label{Stet. Mod.}
     \lim_{\rho \downarrow 0} \omega(\rho) \log\left( \tfrac{1}{\rho}\right) = 0.
\end{equation}

Solutions of systems of the type \eqref{Gleichung} with $p(x)$ growth can typically be shown to be an element of the Sobolev space $W^{m,p(\cdot)}_{loc}(\Omega;\bR^n)$. See Definition \ref{gen.Sob} for an introduction to these spaces.

There have been many investigations on properties of such generalized function spaces in the last years. See for example \cite{Ruzicka:2000}, \cite{Edmunds:2000}, \cite{Harjulehto:2004}, \cite{Diening:2004}, \cite{Haestoe:2006}, and especially \cite{Cruz-Uribe:2003} and \cite{Pick:2001} for properties of the maximal function on generalized Lebesgue spaces. We note that the linear counterpart to this paper, namely the generalization of the classical Calder\'on-Zygmund Theorem \cite{Calderon:1956} to variable Lebesgue spaces has been done by Diening an \Ruzicka\ in \cite{Diening:2003}.

We show in this paper that there exists $\delta > 0$ such that if $|F|^{p(\cdot)}
\in L^q_{loc}(\Omega)$ with $1 < q < \tfrac{n}{n-2}+\delta$, and $u$ is a solution of system \eqref{Gleichung}, then $|D^m u|^{p(\cdot)} \in L^q_{loc}(\Omega)$.

In the case of second order equations ($N = 1$, $m = 1$) and for second order systems with special structure, as for example the $p(x)$ Laplace system, such a result is 
proved in \cite{Acerbi:2004}, without any restriction on $q$. This is due to the fact 
that in this special situation one can prove a $L^{\infty}$ estimate for the 
derivative $Dw$ of the solution $w$ of a suitable frozen problem. In the case of 
general systems (for the second order case see \cite{Kristensen:2006}) this is
not possible. Nevertheless one obtains higher differentiability in the sense
that $D^{m+1} w \in L^{\tilde{p}}$ with a suitable exponent $\tilde{p}$, depending 
on the exponent function $p$. This can be exploited to achieve the desired higher
integrability in the sense of the above statement, with a restriction on the
higher integrability exponent $q$.

The strategy of the proof in this paper follows in a certain sense the ideas in \cite{Acerbi:2004} and \cite{Kristensen:2006}. The key to the proof is an application of a Calder\'on--Zygmund type estimate on level sets of the maximal function of $|D^m u|^{p(\cdot)}$. Therefore the solution will be compared to the solution $w$ of a problem, which is 'frozen' in a point $x_M$ and therefore has the structure of a problem with constant growth exponent $p_2$. The solution $w$ turns out to be higher
differentiable, which translates via Sobolev--Poincar\'e's inequality into higher
integrability of $|D^m w|^{p_2}$. By a suitable comparison estimate between $D^m u$ and $D^m w$, this carries over to the solution $w$.

One should note that we consider both the non degenerate ($\mu \not= 0$) and the degenerate ($\mu = 0$) elliptic case in this paper. Therefore the a priori estimates for the solution of the frozen problem are shown more or less in detail, especially pointing out the differences between the non degenerate and the degenerate case. 

The author should mention that some parts of the proof (especially the comparison
estimate) are widely similar to the proof in the second order case. Therefore at
those points the estimates are shortened very much or cited from other papers. All of the statements are proved in a careful and extensive way in \cite{Habermann:Diss}.

\subsection*{Acknowledgements}

I would like to thank Prof. Dr. Giuseppe Mingione for his helpful comments and many
fruitful discussions about systems and functionals with $p(x)$ growth.

\section{Notations and Setting}\label{sec:setting}

We consider weak solutions of system \eqref{Gleichung}. Note that, using multi indices, \eqref{Gleichung} reads as follows:
\begin{equation*}
     \int_{\Omega} \sum\limits_{|\alpha| = m} A^i_{\alpha} (x,
     D^m u) D^{\alpha} \phi_i \, dx = \int_{\Omega} \sum\limits_{|\alpha| = m} 
     |F(x)|^{p(x)-2}F_{\alpha}^i(x) D^{\alpha}\phi_i\, dx, \quad i=1,\ldots,N.
\end{equation*}
Since the space $\odot^m(\bR^n,\bR^N)$ of symmetric $m$--linear mappings from $\bR^n$ to $\bR^N$ can be identified with the space $\bR^{N\binom{n+m-1}{m}}$, we
consider $A$ as a mapping $A: \Omega \times \bR^{N\binom{n+m-1}{m}} \to \Hom(\bR^{N
\binom{n+m-1}{m}},\bR)$. Additionally for the seek of brevity we introduce the
abbreviation $\cN \equiv N\binom{n+m-1}{m}$. Thus we have $D^m u(x) \in \bR^{\cN}$.

In the sequel we assume that the following structure conditions are satisfied:
Concerning the vector field $A$ we suppose the mapping $z \mapsto A(\cdot,z)$ to be 
of class
$C^0(\bR^{\cN}) \cap C^1(\bR^{\cN})
\setminus \{0\})$ and to satisfy the following growth, ellipticity and
continuity assumptions:
\begin{eqnarray}
     &&\nu \left(\mu^2 + |z|^2\right)^{\frac{p(x)-2}{2}}|\lambda|^2 \le
          \left< D_z A\left(x,z\right) \lambda, \lambda\right> \le L\left(\mu^2+ |z|^2
          \right)^{\frac{p(x)-2}{2}}|\lambda|^2, \label{Wachst. A}\\
     &&\left|A\left(x,z\right) - A\left(y,z\right)\right| \le L \omega\left(
          \left|x-y\right|\right) 
          \Bigl[\left(\mu^2 + |z|^2\right)^{\frac{p(x)-1}{2}} + \left(\mu^2 + |z|^2\right)^{\frac{p(y)-1}{2}}\Bigr]\left|\log \left(\mu^2 + |z|^2\right)\right|, \label{Stetigk. A}
\end{eqnarray}
for all $x,y \in \Omega,\ z,\lambda \in \bR^{\cN}, z \not=0$ where $\nu^{-1},L \in 
\left[1,\infty\right),\ \mu \in \left[0,1\right]$. The parameter $\mu$ is introduced
in order to consider both, the degenerate and the non degenerate case. We assume that
the modulus of continuity $\omega: \bR^+ \to \bR^+$ is a non decreasing, concave and
continuous function satisfying $\omega(0) = 0$.
For the function $p: \Omega \to \left(1,\infty\right)$ we assume that
\begin{equation}\label{Beschr. p}
     1 < \gamma_1 \le p(x) \le \gamma_2 < +\infty,
\end{equation}
for all $x \in \Omega$, as well as
\begin{equation}\label{Stetigk. p}
     \left|p(x) - p(y)\right| \le \omega\left(\left|x-y\right|\right),
\end{equation}
for all $x,y \in \Omega$, where $\omega$ is supposed to fulfill condition \eqref{Stet. Mod.}.

\begin{Bem} By (\ref{Wachst. A}) we can assume that -- eventually enlarging the
constant $L$, reducing $\nu$ respectively -- there holds:
\begin{equation}\label{Wachst. A1}
     \left|A\left(x,z\right)\right| \le L\left(\mu^2+|z|^2\right)^{\left(p(x)-1\right)/2},
\end{equation}
and
\begin{equation}\label{Wachst. A2}
      \nu \left(\mu^2 + |z|^2\right)^{p(x)/2} - L \le \left< A(x,z),z\right> \ \ \ 
      \mbox{ for all }\ x \in \Omega, z \in \bR^{\cN}.
\end{equation}
\hfill $\blacksquare$
\end{Bem}

\begin{definition}[Generalized Lebesgue and Sobolev spaces]\label{gen.Sob}
For a bounded domain $\Omega \subset \bR^n$ and a measurable function $p: \Omega \to
(1,\infty)$ we define the generalized Lebesgue space
\begin{equation*}
L^{p(\cdot)}(\Omega;\bR^N) \equiv \left\{ f \in L^1(\Omega;\bR^N):\ \int_{\Omega}
     |\lambda f(x)|^{p(x)}\, dx < \infty \ \ \mbox{ for some }\ \ \lambda > 0 \right\},
\end{equation*}
which, endowed with the Luxembourg norm
\begin{equation*}
     ||f||_{L^{p(\cdot)}(\Omega;\bR^N)} \equiv \inf\left\{ \lambda > 0: \int_{\Omega} \left|\tfrac{f(x)}{\lambda}\right|^{p(x)}\, dx \le 1\right\}
\end{equation*}
becomes a Banach space. Furthermore the generalized Sobolev space is defined as
\begin{equation*}
     W^{m,p(\cdot)}(\Omega) \equiv \left\{f \in L^{p(\cdot)}(\Omega):\ D^{\alpha} f \in L^{p(\cdot)}(
     \Omega) \quad \mbox{for all}\quad 0 \le |\alpha| \le m \right\},
\end{equation*}
and also becomes a Banach space if endowed with the norm
\begin{equation*}
     ||f||_{W^{m,p(\cdot)}(\Omega;\bR^N)} \equiv \sum_{|\alpha| \le m} ||D^\alpha f||_{L^{p(\cdot)}(\Omega;\bR^N)}.
\end{equation*}
See for example \cite{Ruzicka:2000}, \cite{Harjulehto:2004}, \cite{Cruz-Uribe:2003} and \cite{Edmunds:2000} for more details and further references on these spaces.
\hfill $\blacksquare$
\end{definition}

The main statement of this paper is the following 
\begin{theorem}\label{Hauptsatz 1}
Let $u \in W^{m,p(\cdot)}(\Omega;\bR^N)$ be a weak solution of system (\ref{Gleichung})
under the growth, ellipticity and continuity assumptions (\ref{Wachst. A}) and (\ref{Stetigk. A}) for the vector field $A$, condition (\ref{Beschr. p}) 
for the function $p$ and the assumption (\ref{Stet. Mod.}) for the modulus of
continuity of $p$. Then there exists $\delta \equiv \delta(n,m,\gamma_1,\gamma_2,L/\nu)> 0$ such that if
$|F|^{p(\cdot)} \in L^q_{loc} (\Omega)$ for some
exponent $q$, satisfying
\begin{equation}\label{Bed. q}
1 < q < \frac{n}{n-2} + \delta,
\end{equation}
then
\begin{equation}
     \left|D^m u\right|^{p(\cdot)} \in L^q_{loc}(\Omega).
\end{equation}
In particular there holds: If $\Omega' \Subset \Omega$
and $|F|^{p(\cdot)} \in L^q\left(\Omega'\right)$, then for every given 
$\varepsilon \in \left(0,q-1\right)$ there exists a positive radius $R_0 > 0$, 
depending on
\begin{equation*}
     n,N,m,\gamma_1,\gamma_2,\nu,L,\varepsilon,q,\omega(\cdot),\bigl\| \left|D^m u\right|^{
     p(\cdot)}\bigr\|_{L^1(\Omega)},\bigl\| |F|^{p(\cdot)} \bigr\|_{L^q(\Omega')}
\end{equation*}
such that for any cube $Q_{4R} \Subset \Omega'$ and $R \le R_0$ there holds
\begin{equation*}
     \left( \midint_{Q_R} \left|D^m u\right|^{p(x)q}\, dx\right)^{1/q} \le cK^{ 
     \varepsilon} \midint_{Q_{4R}} \left|D^m u\right|^{p(x)}\, dx + cK^{\varepsilon}
     \left( \midint_{Q_{4R}} |F|^{p(x)q}\, dx + 1\right)^{1/q},
\end{equation*}
where $c \equiv c\left(n,N,m,\gamma_1,\gamma_2,\nu,L,q\right)$ and
\begin{equation*}
     K := \int_{Q_{4R}} \left|D^m u\right|^{p(x)} + |F|^{p(x)(1+\varepsilon)}\, dx +1.
\end{equation*}
\strut\hfill $\blacksquare$
\end{theorem}

\subsubsection*{Remarks on the notation}

In the whole paper
$\Omega \subset \bR^n, (n\ge 2)$ denotes a bounded domain in the space $
\bR^n$ and $Q\left(x,R\right) \equiv Q_R(x)$ a cube whose axes are parallel to the
axes of the coordinate system, with center $x$ and side length $2R$. Sometimes
$R$ will also be called the 'radius' of the cube. The Lebesgue measure of a measurable 
set $A$ is abbreviated by $|A| \equiv {\mathcal L}^n(A)$. For a locally integrable function $u \in L^1_{loc}(\Omega)$ we define the mean value on the cube $Q$ by
\begin{equation*}
     \left(u\right)_{x_0,R} := \midint_{Q\left(x_0,R\right)} 
        u(x)\, dx = \frac{1}{|Q(x_0,R)|}\int_{Q\left(x_0,R\right)} 
        u(x)\, dx.
\end{equation*}
In the case the centre of the cube is obvious from the context, we will often
just write $Q_R$ or $Q$ instead of $Q(x_0,R)$, $(u)_R$ instead of $(u)_{x_0,R}$
respectively.

The letter $c$ denotes a constant which will not necessarily be the same at different places
in the paper and which may sometimes change from line to line. Constants that will
be referred to at other points of the work, will be signed in a unique way, mostly by
different indices. In the case we want to emphasize the fact that a constant changes
from one line to another, we will label this by mathematical accents, as for example
$\tilde{c}$ or $\bar{c}$. For the survey we will not specify the dependencies of the 
constants in between the estimates, but of course at the end of them.

For $\Omega \subset \bR^n, p > 1$, let $L^{p}(\Omega;\bR^N)$ be the well known Lebesgue
space to the power $p$. For $m \in \bN$ we define the Sobolev space 
\begin{equation*}
     W^{m,p}\left( \Omega; \bR^N \right) := \left\{ u \in L^p\left( \Omega, 
     \bR^N \right): D^{\alpha}u \in L^p\left( \Omega\right) \mbox{ for } 0 
     \le \left|\alpha\right|\le m \right\},
\end{equation*}
with the multi-index $\alpha = (\alpha_1,\ldots,\alpha_n) \in \bN^{n}$ and the
abbreviations $|\alpha| := \alpha_1 + \ldots + \alpha_n$ and $D^{\alpha} u :=
D^{\alpha_1}_1\ldots D^{\alpha_n}_n u$. Furthermore let $W^{m,p}_0(\Omega;\bR^N)$
denote the closure of $C^{\infty}(\Omega;\bR^N)$ in the space $W^{m,p}(\Omega;\bR^N)$.

\section{A priori estimates}

As we will see in the proof, the most important difficulty compared to the proof in the second
order case ($m = 1$, see \cite{Acerbi:2004}) comes from the a priori estimate, which 
is different for higher order systems, since one can not obtain an $L^{\infty}$-- bound
for the derivative $D^m u$. The optimal result is the following
\begin{Lem}\label{A priori Lemma}
Let $\Omega \subset \bR^n$ be a bounded domain, $p > 1$ constant and 
$w \in W^{m,p} (\Omega; \bR^N)$ a weak solution of the system
\begin{equation}\label{Gl.Konst.}
     \int_{\Omega} \left<A\left(D^m w\right), D^m \phi\right>\, dx = 0 
      \quad \quad \mbox{ for all }\ \phi \in W_0^{m,p}\left(\Omega;\bR^N\right),
\end{equation}
in which the function $A: \bR^{\cN} \to \Hom(\bR^{\cN},\bR)$ is of the class $C^1$ and satisfies the conditions
\begin{equation}\label{Bed.A}
\begin{aligned}
     \nu\left(\mu^2 + |z|^2\right)^{\frac{p-2}{2}}|\lambda|^2 &\le \left< DA(z)\lambda,\lambda\right>
     \le L\left(\mu^2 + |z|^2\right)^{\frac{p-2}{2}}|\lambda|^2,\\
     |A(z)| &\le L\left(\mu^2 + |z|^2\right)^{\frac{p-1}{2}},
\end{aligned}
\end{equation}
for all $z\in \bR^{\cN}$. Then the following holds:\\[0.1cm]
In the case $1 < p < 2$ we have $w \in W^{m+1,p}_{loc}\left(\Omega;\bR^N\right)$ 
together with the estimate
\begin{equation}\label{Apr.p<2.1}
     \int_{Q_{\tau R}} \left|D^{m+1}w\right|^p\, dx 
     \le \frac{c}{R^2} \int_{Q_R} \left( 
           \mu^2 + \left|D^m w\right|^2\right)^{p/2}\, dx.
\end{equation}
Additionally distinguishing the cases $\mu \not= 0$ and $\mu = 0$, we obtain
furthermore
\begin{eqnarray}
     \int_{Q_{\tau R}} \left| D\left[\left(\mu^2 + |D^m w|^2\right)^{p/4}\right] 
     \right|^2\, dx
     \le \frac{c}{R^2} \int_{Q_R} \left(\mu^2 + |D^m w|^2\right)^{p/2}\, dx
     && (\mu \not= 0) \label{Apr.p<2.2}\\
     \int_{Q_{\tau R}} \left| D\left[ |D^m w|^{\frac{p-2}{2}}D^m w\right]\right|^2
     \, dx \le \frac{c}{R^2} \int_{Q_R} |D^m w|^p\, dx
     && (\mu = 0).\label{Apr.p<2.3}
\end{eqnarray}
for any cube $Q_R \Subset \Omega$ and any $\tau \in (0,1)$.\\[0.1cm]
In the case $p \ge 2$
we have $D\bigl[\left(\mu^2 + |D^m w|^2\right)^{\frac{p-2}{4}}D^m w\bigr] \in L^{2}_{loc}
(\Omega)$ and for any cube $Q_R \Subset \Omega$ and any $\tau \in (0,1)$ there holds
\begin{equation}\label{Apr.p>2.1}
     \int_{Q_{\tau R}} \left| D\left[\left(\mu^2 + |D^m w|^2\right)^{\frac{p-2}{4}}D^m 
     w\right]\right|^2
     \, dx \le \frac{c}{R^2} \int_{Q_R} \left( \mu^2 
     + |D^m w|^2\right)^{p/2}\, dx.
\end{equation}
Moreover the constants in the estimates above depend on $n,N,m,p,\tau$ and $L/\nu$ and
$c \uparrow \infty$ as $\tau \uparrow 1$.
\strut\hfill $\blacksquare$
\end{Lem}

\begin{proof}
We start by proceeding analoguously to the second order case. For $h \in \bR$ with
$|h| < (1-\tau)R$ and $x \in Q_{\tau R}$  we denote by
\begin{equation*}
     \tau_{s,h}w(x) := w\left(x+he_s\right) - w\left(x\right)
\end{equation*}
the finite difference and by
\begin{equation*}
     \Delta_{s,h} w := \frac{\tau_{s,h} w}{h},
\end{equation*}
the difference quotient in direction $e_s$, where 
$e_s$ denotes the $s$th unit vector in $\bR^n$. We consider
the test function
\begin{equation*}
     \phi = \Delta_{s,-h} \left(\eta^{2m} \Delta_{s,h} (w-P)\right),
\end{equation*}
with a suitable polynomial $P$.
Moreover for $R > 0$ and $0 < |h| < R$ we denote by
\begin{equation*}
     Q_{R,h} \equiv \left\{ x \in Q_R:\ \mbox{dist} \left(x, \partial Q_R 
     \right) > |h| \right\}
\end{equation*}
the inner parallel cube, whose sidelength is $R - |h|$.
For $w \in W^{m,p}\left(Q_R\right)$ we have $\Delta_{s,h} w \in 
W^{m,p}\left(Q_{R,h}\right)$.
In the definition of $\phi$ we choose $\eta \in C^{\infty}_c \left(Q_R\right)$ to be a standard cut-off 
function with the properties 
\begin{equation*}
     0 \le \eta \le 1,\ \ \eta \equiv 1 \ \mbox{ on } Q_{\tau R},\ \ \mbox{spt} 
     \eta \Subset Q_{\sqrt{\tau}R}
\end{equation*}
and
\begin{equation*}
     \left| D^k \eta \right| \le \frac{c}{(R(\sqrt{\tau}-\tau))^k} = 
     \frac{c(\tau)}{R^k}, \ \ \ \mbox{ for } \ \ k = 1,\ldots,m,
\end{equation*}
with $c(\tau) \to \infty$ as $\tau \nearrow 1$ or $\tau \searrow 0$.
By the chain rule we immediately get
\begin{equation}\label{Abschn.2}
     \left|D^k \left(\eta^{2m}\right)\right| \le c(n,m,\tau) R^{-k} \sum\limits_{
     j=1}^k \eta^{2m-j} \le c(n,m,k,\tau) R^{-k} \eta^{2m-k}.
\end{equation}
$P:\Omega \to \bR^N$ denotes the unique polynomial of degree $m-1$ whose coefficients
are chosen to satisfy
\begin{equation}\label{MW Polynom}
     \left(D^k\left(w-P\right)\right)_{Q_{\sqrt{\tau}R}} = \midint_{Q_{\sqrt{\tau}R}} 
     \mskip-15mu D^k\left(w-P\right)\, dx = 0 \ \ \mbox{ for }\ k=0,\ldots,m-1.
\end{equation}
Existence and uniqueness of such polynomials are well known and can be found for example in \cite{Giaquinta:1979}.
Testing \eqref{Gl.Konst.}, using standard identities for difference quotients and the
general chain rule, we obtain for $|h| \le R(1-\sqrt{\tau})$:
\begin{eqnarray*}
     &&0  = \int_{Q_{\sqrt{\tau}R}} \left< \Delta_{s,h}A, \eta^{2m}D^m\left(
          \Delta_{s,h}w\right)\right>\, dx \\
     &&\mskip+70mu +
          \int_{Q_{\sqrt{\tau}R}} \Bigl<\Delta_{s,h}A, \sum\limits_{k=1}^{m} 
	  \binom{m}{k}D^k\left(\eta^{2m}\right)\odot D^{m-k}\left(\Delta_{s,h}\left(w-P
          \right)\right)\Bigr>\, dx\\
     &&\mskip+14mu = I^{(1)} + I^{(2)}.
\end{eqnarray*}
Now we distinguish the cases $1 < p < 2$ and $p \ge 2$. We start with the 
case $p \ge 2$. We translate the growth and ellipticity conditions for $A$ into conditions for $\Delta_{s,h} A$. Elementary calculations together with the differentiability of $A$ show that 
\begin{equation}\label{Ident.B_h}
     \Delta_{s,h}A(x) =  B_h(x) D^m \Delta_{s,h} w(x),
\end{equation}
with
\begin{equation*}
     B_h \equiv B_h(x) := \int_0^1 DA\left(D^m w(x) + t D^m \tau_{s,h}w(x) 
          \right)\, dt.
\end{equation*}
Using \eqref{Bed.A} and the technical lemma \ref{techn. Lemma} (with
exponent $\frac{p-2}{2}$) we obtain the following pointwise 
estimates for $B_h$:
\begin{equation}\label{Eigensch.B}
     \left|B_h\right| \le c_1(p) L W_h^{p-2},\qquad \left<B_h\lambda, \lambda\right> 
     \ge c_2(p) \nu W_h^{p-2} \left|\lambda\right|^2,
\end{equation}
for $\lambda \in \bR^{\cN}$, with $W_h^2 := \mu^2 + |D^m w|^2 + |\tau_{s,h} D^m w|^2$.

Now, $I^{(1)}$ ist estimated by $(\mbox{\ref{Eigensch.B}})_2$ from below as follows:
\begin{equation*}
     I^{(1)} \ge c_2 \nu \int_{Q_{\sqrt{\tau}R}} \eta^{2m} W_h^{p-2} \left| \Delta_{s,h} 
           D^m w \right|^2\, dx.
\end{equation*}
From (\ref{Ident.B_h}) and (\ref{Abschn.2}) we infer that 
\begin{equation*}
     I^{(2)} 
      \le  c \sum\limits_{k=1}^{m} \int_{Q_{\sqrt{\tau}R}} \left|B_h\right| \left|
           \Delta_{s,h} D^m w\right| R^{-k}\eta^{2m-k}\left|\Delta_{s,h} 
           D^{m-k} \left(w-P\right)\right|\, dx
      =:   c \sum_{k=1}^m I^{(2)}_k,
\end{equation*} 
where $c \equiv c(n,m,\tau)$. $(\mbox{\ref{Eigensch.B}})_1$ and Young's inequality
lead to
\begin{eqnarray*}
     I^{(2)}_k 
     &\le& c_1 L \varepsilon \int_{Q_{\sqrt{\tau}R}} \eta^{2m} W_h^{p-2} \left|
          \Delta_{s,h} D^m w\right|^2 \, dx \\
     &   &\mskip+40mu + \frac{c_1L}{4\varepsilon} \int_{Q_{\sqrt{\tau}R}}R^{-2k}
          \eta^{2\left( m-k\right)}W_h^{p-2} \left|\Delta_{s,h} D^{m-k}\left( 
          w-P\right) \right|^2 \, dx.
\end{eqnarray*} 
Summing up the estimates for $I^{(2)}_k$, subsequently choosing $\varepsilon = \frac{c_2 \nu}{2c_3 Lm}$ and taking into account that $\eta \equiv 1$ on $Q_{\tau R}$ finally leads to
\begin{equation}\label{Absch.1}
     \int_{Q_{\tau R}} W_h^{p-2} \left|\Delta_{s,h} D^m w\right|^2\, dx
     \le c\sum\limits_{k=1}^m R^{-2k} \int_{Q_{\sqrt{\tau}R}} W_h^{p-2} 
          \left|\Delta_{s,h} D^{m-k}\left(w-P\right)\right|^2\, dx,
\end{equation} 
where the constant $c$ depends on $n,m,p,\tau$ and $L/\nu$. Now we are going to estimate each of the terms appearing on the right hand side of
(\ref{Absch.1}). Writing
$R^{-2k} = R^{-2\left(p-2\right)/p} \cdot R^{\left(2(1-k)p-4\right)/p}$, we obtain by
Young's inequality (with exponents $p \equiv \frac{p}{p-2} > 1,\ q \equiv \frac{p}{2}$) for every $k=1,\ldots,m$
\begin{eqnarray*}
     &&R^{-2k}\int_{Q_{\sqrt{\tau}R}}\mskip-10mu W_h^{p-2} \left|\Delta_{s,h} 
          D^{m-k}\left(w-P\right)\right|^2\, dx\\
     &&\mskip+65mu \le \tfrac{p-2}{p}R^{-2}\int_{Q_{\sqrt{\tau}R}}\mskip-10mu 
          W_h^{p}\, dx + \tfrac{2}{p} R^{p(1-k)-2} \int_{Q_{\sqrt{\tau}R}}\mskip-10mu 
          \left| \Delta_{s,h} D^{m-k}\left(w-P\right)\right|^{p}\, dx.
\end{eqnarray*}
Since $w \in W^{m,p}\left(Q_R\right)$, using standard estimates for difference quotients (note that $|h| \le R(1-\sqrt{\tau})$) and subsequently applying Poincar\'e's inequality $(k-1)$ times (note the choice of the polynomial $P$ in (\ref{MW Polynom})) provides for any $k=1,\ldots,m$:
\begin{equation*}
     \int_{Q_{\sqrt{\tau}R}} \left|\Delta_{s,h}D^{m-k}\left(w-P\right)\right|^{p}
          \, dx
     \le  \int_{Q_{R}} \left|D^{m-k+1}\left(w-P\right) 
          \right|^{p}\, dx
     \le  cR^{(k-1)p} \int_{Q_{R}} \left|D^m w\right|^{p}\, dx,
\end{equation*}
which finally leads to
\begin{equation*}
     \int_{Q_{\tau R}} W_h^{p-2}\left|\Delta_{s,h} D^m w\right|^2\, dx \le c \sum
     \limits_{k=1}^m \Bigl[ R^{-2} \int_{Q_{\sqrt{\tau}R}} W_h^{p}\, dx
     + R^{-2} \int_{Q_{R}} \left|D^m w\right|^{p}\, dx \Bigr],
\end{equation*}
where the constant $c$ depends only on $n,N,m,p,\tau$ and $L/\nu$.
Furthermore we easily see
\begin{equation*}
     \int_{Q_{\sqrt{\tau}R}} W_h^{p}\, dx 
     \le c(p) \int_{Q_{R}} \left(\mu^2+\left|D^m w\right|^2\right)^{
         p/2}\, dx.
\end{equation*}
Therefore we obtain for any $s=1,\ldots,n$ and $0 < |h| \le R(1-\sqrt{\tau})$ that
\begin{equation}\label{Absch.2}
     \int_{Q_{\tau R}} W_h^{p-2}\left|\Delta_{s,h}D^m w\right|^2\, dx \le
     cR^{-2} \int_{Q_{R}} \left(1+\left|D^m w\right|^2\right)^{p/2}\, dx.
\end{equation}
with $c \equiv c(n,N,m,p,\tau,L/\nu)$.
In the case $1 < p < 2$ we proceed in a different way. However we note that the arguments which lead to the bound from below for $I^{(1)}$ in the
case $p \ge 2$ also work here. Using $(\mbox{\ref{Eigensch.B}})_2$ we get
\begin{equation*}
     I^{(1)} \ge c_2\nu \int_{Q_{\sqrt{\tau}R}} \eta^{2m}W_h^{p-2}\left|\Delta_{s,h} D^m w
     \right|^2\, dx.
\end{equation*}
To treat $I^{(2)}$ we use a different formula for $\Delta_{s,h} A$. We do the following
formal calculation, which holds for functions $w \in W^{m+1,p}$. The result for
$w \in W^{m,p}$ can then be achieved by approximation. We write
\begin{equation}\label{Ident.Bt_h}
     \left[\Delta_{s,h}A(D^m w)\right](x) 
     = \frac{1}{h}\int_0^1 \frac{d}{dt} A\left(D^m w\left(x+the_s
          \right)\right)\, dt
     = D_s \tilde{B}_h(x),
\end{equation}
with
\begin{equation*}
     \tilde{B}_h(x) \equiv \int_0^1 A\left(D^m w\left(x+the_s\right)\right)\, dt.
\end{equation*}
By the growth condition (\ref{Bed.A}) for $A$ we find that
\begin{equation*}
     |\tilde{B}_h| \le L\int_0^1\left(\mu^2+\left|D^m u\left(x+the_s
     \right)\right|^2\right)^{\frac{p-1}{2}} =: L \cdot Y_h.
\end{equation*}
For $I^{(2)}$ we write
\begin{equation*}
     I^{(2)} = c(m)\sum\limits_{k=1}^m I^{(2)}_k,
\end{equation*}
with
\begin{equation*}
     I^{(2)}_k \equiv \int_{Q_{\sqrt{\tau}R}} \left<D_s \tilde{B}_h,D^k\left(\eta^{2m} 
     \right) \odot D^{m-k}\left(\Delta_{s,h}\left(w-P\right)\right)\right>
     \, dx.
\end{equation*}
Taking into account (\ref{Ident.Bt_h}) and $\mbox{spt} \eta \Subset Q_{\sqrt{\tau}R}$ we
obtain for $I^{(2)}_k$ by partial integration 
\begin{eqnarray*}
     &&I^{(2)}_k  = - \int_{Q_{\sqrt{\tau}R}} \left<\tilde{B}_h, D_s\left(D^k\left( 
          \eta^{2m}\right) \right)D^{m-k}\left(\Delta_{s,h}\left(w-P\right)\right)
          \right>\, dx\\
     &&\mskip+145mu - \int_{Q_{\sqrt{\tau}R}} \left<\tilde{B}_h, D^k\left(\eta^{2m} 
          \right) D_s D^{m-k}\left(\Delta_{s,h}\left(w-P\right)\right)\right>\, 
          dx.
\end{eqnarray*}
By (\ref{Abschn.2}) we obtain
\begin{equation*}
     \left|D_sD^k\left(\eta^{2m}\right)\right| \le \left|D^{k+1}\left(
          \eta^{2m}\right)\right| \le c(n,m,k,\tau) R^{-(k+1)} \eta^{2m-(k+1)},
\end{equation*}
and therefore
\begin{eqnarray*}
     &&|I^{(2)}_k| \le cR^{-(k+1)} \int_{Q_{\sqrt{\tau}R}} Y_h \eta^{2m-(k+1)} \left|
          D^{m-k}\Delta_{s,h} \left(w-P\right)\right|\, dx\\
     &&\mskip+100mu + cR^{-k} \int_{Q_{\sqrt{\tau}R}} Y_h \eta^{2m-k}\left|D^{m-k+1} 
          \Delta_{s,h} \left(w-P\right)\right|\, dx,
\end{eqnarray*}
with constants $c \equiv c(n,m,\tau)$. Combining the previous estimates we arrive at:
\begin{equation}\label{Absch.3}
\begin{aligned}
     c_2\nu\int_{Q_{\sqrt{\tau}R}} W_h^{p-2}\eta^{2m} \left|\Delta_{s,h}D^m 
          w\right|^2\, dx 
     &\le c\sum\limits_{k=1}^m R^{-(k+1)} \int_{Q_{\sqrt{\tau}R}} Y_h 
          \eta^{2m-(k+1)}\left|D^{m-k}\Delta_{s,h} \left(w-P\right)\right|\, dx\\
     &\mskip+15mu + c\sum\limits_{k=1}^m R^{-k} \int_{Q_{\sqrt{\tau}R}} Y_h 
          \eta^{2m-k}\left|D^{m-k+1}\Delta_{s,h}\left(w-P\right)\right|\, dx\\
     &= cR^{-1}\int_{Q_{\sqrt{\tau}R}} Y_h \eta^{2m-1} \left|D^m \Delta_{
          s,h} w\right|\, dx\nonumber\\
     &\mskip+15mu + c\sum\limits_{k=1}^m R^{-(k+1)} \int_{Q_{\sqrt{\tau}R}} Y_h 
          \eta^{2m-(k+1)} \left|D^{m-k}\Delta_{s,h} \left(w-P\right)\right|\, dx\\
     &= I^{(3)} + \sum\limits_{k=1}^m I^{(4)}_k.
\end{aligned}
\end{equation}
We first consider $I^{(3)}$. Applying Young's inequality we get
\begin{eqnarray*}
     R^{-1} Y_h \eta^{2m-1}\left|D^m \Delta_{s,h} u\right|
     & = &R^{-1}Y_hW_h^{\frac{2-p}{2}}W_h^{\frac{p-2}{2}}\eta^m
          \eta^{m-1}\left|D^m \Delta_{s,h} w\right|\\[0.1cm]
     &\le&\varepsilon \eta^{2m}W_h^{p-2}\left|D^m \Delta_{s,h} w
          \right|^2 + \frac{1}{4\varepsilon} R^{-2} \eta^{2m-2} W_h^{2-p}Y_h^2.
\end{eqnarray*}
By a suitable choice of $\varepsilon$, we can absorb the first term on the
left hand side of (\ref{Absch.3}). The second term can be estimated by Young's
inequality (with exponents $p \equiv \frac{p}{2-p}$, $q \equiv \frac{p}{2 
\left(p-1\right)}$) as follows:
\begin{equation*}
     R^{-2} W_h^{2-p}Y_h^2 \le R^{-2}\left(\tfrac{2-p}{p}\ W_h^{p} + 
     \tfrac{2(p-1)}{p}\ Y_h^{\frac{p}{p-1}}\right) \le c(p)R^{-2}\left(W_h^{p} + 
     Y_h^{\frac{p}{p-1}}\right).
\end{equation*}
To estimate $I_k^{(4)}$ we write $R^{-(k+1)} = R^{2\left(1-p\right)/p}
R^{\left(p(1-k)-2\right)/p}$ and obtain by Young's inequality
\begin{eqnarray*}
     I^{(4)}_k 
     & = &cR^{-(k+1)}\int_{Q_{\sqrt{\tau}R}}\eta^{2m-(k+1)} Y_h \left|D^{m-k} 
          \Delta_{s,h}\left(w-P\right)\right|\, dx\\
     &\le&c\tfrac{p-1}{p} R^{-2} \int_{Q_{\sqrt{\tau}R}} \eta^{2m-(k+1)} 
          Y_h^{\frac{p}{p-1}}\, dx\\
     &&\mskip+40mu +c\tfrac{1}{p}R^{p(1-k)-2} \int_{Q_{\sqrt{\tau}R}}\eta^{2m-(k+1)} 
          \left|D^{m-k} \Delta_{s,h} \left(w-P\right)\right|^{p}\, dx\\
     &\le&cR^{-2}\int_{Q_{\sqrt{\tau}R}}\mskip-10mu Y_h^{\frac{p}{ 
          p-1}}\, dx + cR^{p(1-k)-2}\int_{Q_{\sqrt{\tau}R}} \left|D^{m-k} 
          \Delta_{s,h} \left(w-P\right) \right|^{p}\, dx.
\end{eqnarray*}
Combining the previous estimates we arrive at
\begin{eqnarray*}
     \int_{Q_{\tau R}} W_h^{p-2}\eta^{2m} \left|\Delta_{s,h}D^m 
          w\right|^2\, dx
     &\le&cR^{-2} \int_{Q_{\sqrt{\tau}R}}\mskip-10mu 
          W_h^{p}\, dx + cR^{-2}\int_{Q_{\sqrt{\tau}R}} \mskip-10mu 
          Y_h^{\frac{p}{p-1}}\, dx\\
     &&   + c\sum\limits_{k=1}^m 
          R^{p(1-k)-2}\int_{Q_{\sqrt{\tau}R}}\mskip-10mu \left|D^{m-k} 
          \Delta_{s,h}\left(w-P\right)\right|^{p}\, dx,
\end{eqnarray*}
with constants $c \equiv c(n,m,\tau,p,L/\nu)$.
By H\"older's inequality and Fubini's theorem we see
\begin{equation*}
     \int_{Q_{\sqrt{\tau}R}} Y_h^{\frac{p}{p-1}}\, dx 
          \le \int_{Q_{\sqrt{\tau}R}} \int_0^1\left(\mu^2+\left|D^m w\left(x + 
          the_s\right)\right|^2\right)^{p/2}\, dt\, dx
          \le \int_{Q_{R}} \left(\mu^2+\left|D^m w(x)\right|^2 
          \right)^{p/2}\, dx.
\end{equation*}
As in the case $p \ge 2$ we deduce easily
\begin{equation*}
     \int_{Q_{\sqrt{\tau}R}} W_h^{p}\, dx \le c(p)\int_{Q_{R}} \left(\mu^2+\left|D^m w
     \right|^2\right)^{p/2}\, dx.
\end{equation*}
Furthermore again by standard estimates for difference quotients
we obtain
\begin{equation*}
     \int_{Q_{\sqrt{\tau}R}} \left|D^{m-k}\Delta_{s,h}\left(w-P\right)\right|^{ 
          p}\,dx 
     \le c(n,N)\int_{Q_R} \left|D^{m-k+1}\left(w-P\right)\right|\, dx.
\end{equation*}
Inserting this above we find that (note that $\eta \equiv 1$ on $Q_{\tau R}$)
\begin{eqnarray*}
     \int_{Q_{\tau R}} W_h^{p-2} \left|\Delta_{s,h}D^m w\right|^2\, dx
     &\le&cR^{-2}\int_{Q_{R}} \left(\mu^2+\left|D^m 
          w\right|^2 \right)^{p/2}\, dx\\
     &&\mskip+50mu + c\sum\limits_{k=1}^m R^{p(1-k)-2}\int_{Q_R} 
          \left|D^{m-k+1} \left(w-P\right) \right|^{p}\, dx,
\end{eqnarray*}
with $c \equiv c(n,N,m,p,\tau,L/\nu)$.
By the choice of the polynomial $P$ (see \eqref{MW Polynom}) we can apply $(k-1)$
times Poincar\'e's inequality to the integrals $\int_{Q_R} |D^{m-k+1}(w-P)|^p\, dx$;
actually we have that 
\begin{equation*}
     \int_{Q_{\sqrt{\tau}R}} \left|D^{m-k+1} \left(w-P\right)\right|^{p}\, dx 
     \le c(n,N) R^{(k-1)p} \int_{Q_{\sqrt{\tau}R}} \left|D^m w \right|^{p}\, dx,
\end{equation*}
for all $k = 2,\ldots,m$.
Therefore we obtain
\begin{eqnarray*}
     \int_{Q_{\tau R}} W_h^{p-2} \left|\Delta_{s,h} D^m w\right|^2\, dx
     &\le&cR^{-2} \int_{Q_{R}} \left(\mu^2+\left|D^m 
          w\right|^2\right)^{p/2}\, dx
     + c\sum\limits_{k=1}^m R^{-2} \int_{Q_{\sqrt{\tau}R}} 
          \left|D^m w\right|^{p}\, dx\\
     &\le&cR^{-2} \int_{Q_{R}} \left(\mu^2+\left| D^m w\right|^2\right)^{p/2}\, dx,
\end{eqnarray*}
with $c\equiv c(n,N,m,\tau,p,L/\nu)$.
Hence for any $p > 1$ there holds
\begin{equation}\label{Absch. 4}
     \int_{Q_{\tau R}} W_h^{p-2} \left|\Delta_{s,h} D^m w\right|^2\, dx \le 
     cR^{-2} \int_{Q_{R}} \left(\mu^2+\left|D^m w\right|^2 \right)^{p/2}\, dx,
\end{equation}
with a constant $c \equiv c(n,N,m,p,\tau,L/\nu)$.

Now we distinguish the cases
$1<p<2$ and $p\ge2$:

\noindent
The case $1 < p < 2$: 
We set $2 \alpha :=p\left(2-p\right)$, obtaining by 
Young's inequality 
\begin{equation*}
     \left|\Delta_{s,h} D^m w\right|^{p} = W_h^{\alpha} W_h^{-\alpha} \left|
     \Delta_{s,h} D^m w\right|^{p} \le c(p)\left(W_h^{p} + W_h^{p-2} \left|
     \Delta_{s,h}D^m w\right|^2\right),
\end{equation*}
and therefore
\begin{equation}\label{Delta D^mw bd.}
\int_{Q_{\tau R}} \left|\Delta_{s,h} D^m w\right|^{p}\, dx \le cR^{-2}
     \int_{Q_{R}} \left(\mu^2+\left|D^m w\right|^2\right)^{p/2}\, dx.
\end{equation}

By \eqref{Delta D^mw bd.}, we see that the 
sequence $\Delta_{s,h} D^m w$ is uniformly bounded in $L^{p}\left( 
Q_{\tau R}\right)$. Therefore $\Delta_{s,h} D^m w$ converges as $h \to 0$ strongly in $L^{p}_{loc}\left(Q_{\tau R}\right)$ to $D_sD^m w$, i.e. $w \in W^{m+1,p}_{loc} \left(Q_{\tau R}\right)$ and
\eqref{Apr.p<2.1} holds. On the other hand, for $\mu \in (0,1]$, a subsequence of
$\Delta_{s,h}D^m w$ converges pointwise a.e. to $D_sD^m w$ as $h \to \infty$. 
With the convergence
\begin{equation*}
      \tau_{s,h} D^m w \overset{h \to 0}{\longrightarrow} 0 \ \ \ \mbox{ in } 
      L^{p} \left(Q_{\tau R}\right),
\end{equation*}
there holds
\begin{equation*}
      W_h = \left(\mu^2+\left|D^m w\right|^2 + \left|\tau_{s,h} D^m w 
      \right|^2\right)^{1/2} \overset{h \to 0}{\longrightarrow} \left(\mu^2+
      \left|D^m w\right|^2\right)^{1/2} \ \ \mbox{ in } L^{p}\left(Q_{\tau R}\right).
\end{equation*}
This implies the pointwise almost everywhere convergence of a subsequence of $W_h$ to
$(\mu^2+ \linebreak |D^m w|^2)^{1/2}$. By Fatou's Lemma we now conclude 
with \eqref{Absch. 4}:
\begin{equation*}
     \int_{Q_{\tau R}} \left(\mu^2 + |D^m w|^2\right)^{\frac{p-2}{2}}|D_sD^m w|^2\, dx
     \le  \frac{c}{R^2} \int_{Q_R} \left(\mu^2 + |D^m w|^2\right)^{p/2}\, dx.
\end{equation*}
Therefore we end up with
\begin{equation}
     \int_{Q_{\tau R}} \left(\mu^2 + |D^m w|^2\right)^{\frac{p-2}{2}}|D^{m+1} w|^2
     \, dx \le \frac{c}{R^2} \int_{Q_R} \left(\mu^2 + |D^m w|^2\right)^{p/2}\, dx.
\end{equation}
On the other hand, by differentiating, we have the estimate
\begin{equation}\label{CZ.Apr.est.diff}
     \left|D\left[\left(\mu^2 + |D^m w|^2\right)^{p/4}\right]\right|^2
     \le c\left(\mu^2 + |D^m w|^2\right)^{\frac{p-2}{2}}|D^{m+1}w|^2.
\end{equation}
Therefore we conclude
\begin{equation*}
     \int_{Q_{\tau R}} \left|D\left[\left(\mu^2 + |D^m w|^2\right)^{p/4}\right]
     \right|^2\, dx \le \frac{c}{R^2} \int_{Q_R} \left(\mu^2 + |D^m w|^2
     \right)^{p/2}\, dx,
\end{equation*}
which is exactly \eqref{Apr.p<2.2}. We use an elementary algebraic property of the function $V_\mu (z) := (\mu^2+|z|^2)^{(p-2)/4}z$ (see \eqref{prop.V} on page \pageref{prop.V}) to obtain 
\begin{equation}\label{CZ.Apr.est.techn.1}
\begin{aligned}
     &\left| \tau_{s,h} \left[ \left( \mu^2 + |D^m w|^2\right)^{\frac{p-2}{4}} D^m w
     \right] \right|^2\\
     &\mskip+50mu
     \le c \left[ \mu^2 + \left|D^m w(x)\right|^2 + \left|D^m w(x+he_s)\right|^2 
     \right]^{\frac{p-2}{2}}\left| D^m w(x+he_s) - D^m w(x)\right|^2\\
     &\mskip+50mu
     \le c \left[ \mu^2 + \left|D^m w\right|^2 + \left|\tau_{s,h}D^m w\right|^2 
     \right]^{\frac{p-2}{2}}\left|\tau_{s,h} D^m w\right|^2.
\end{aligned}
\end{equation}
Combining this estimate for $\mu = 0$ with \eqref{Absch. 4}, we end up with
\begin{equation*}
     \int_{Q_{\tau R}} \left|\Delta_{s,h} \left[ |D^m w|^{\frac{p- 
     2}{2}} D^m w \right] \right|^2 \le \frac{c}{R^2} \int_{Q_R} |D^m 
     w|^{p}\, dx.
\end{equation*}
We see that the sequence $\Delta_{s,h} [ |D^m w|^{\frac{p- 
2}{2}} D^m w ]$ is uniformly bounded in $L^2(Q_{\tau R})$. By a standard lemma 
about difference quotients it converges as $h \to \infty$ strongly in $L^2_{loc}(Q_{
\tau R})$ to $D_s [ |D^m w|^{\frac{p- 2}{2}} \linebreak D^m w ]$. The estimate
above together with the convergence yield the desired estimate \eqref{Apr.p<2.3}.

\noindent
The case $p \ge 2$: 
We take \eqref{CZ.Apr.est.techn.1} together with \eqref{Absch. 4} to conclude
\begin{equation}\label{Delta D^mw bd.2}
    \int_{Q_{\tau R}} \left| \Delta_{s,h} \left[ \left( \mu^2 + |D^m w|^2\right)^{
     \frac{p-2}{4}} D^m w \right] \right|^2 \le cR^{-2} \int_{Q_R} \left( \mu^2 + 
     \left|D^m w\right|^2\right)^{p/2}\, dx.
\end{equation}
By \eqref{Delta D^mw bd.2} we see that the sequence
$\Delta_{s,h} [ ( \mu^2 + |D^m w|^2)^{\frac{p-2}{4}} D^m w ] $ is 
uniformly bounded in $L^2(Q_{\tau R})$ and therefore the sequence converges as $h \to 0$ strongly in $L^2_{loc}(
Q_{\tau R})$ to $D_s [ ( \mu^2 + |D^m w|^2)^{\frac{p-2}{4}} D^m w ]$.
The estimates above together with the convergence provide the desired estimate
\eqref{Apr.p>2.1}.
\end{proof}

\section{An additional Gehring improvement}
Starting by Lemma \ref{A priori Lemma}, we can now achieve by standard techniques a 
further higher integrability exponent in the following sense:

\begin{Lem}\label{A priori Gehring}
Let $w \in W^{m,p}(\Omega;\bR^N)$ be a solution of the system \eqref{Gl.Konst.}, which
satisfies the structure conditions \eqref{Bed.A}. Then there exists $\delta \equiv
\delta(n,m,p,L/\nu)$ and a constant $c \equiv c(n,m,p,L,M)$ such that the following
holds:

In the case $1 < p < 2$ and $\mu \not= 0$ we have the estimate
\begin{equation}\label{HI.Dm+1.1}
     \int_{Q_{\tau R}} \left| D\left[(\mu^2 + |D^m w|^2)^{p/4}\right]\right|^{2(1+
     \delta)}\, dx \le \left( \frac{c}{R^2} \int_{Q_R} (\mu^2 + |D^m 
     w|^2)^{p/2}\, dx \right)^{1+\delta}.
\end{equation}

In the case $1 < p < 2$ and $\mu = 0$ we obtain
\begin{equation}\label{HI.Dm+1.2}
     \int_{Q_{\tau R}} \left| D\left[|D^m w|^{\frac{p-2}{2}}D^m w\right]\right|^{2(
     1+\delta)}\, dx \le \left( \frac{c}{R^2} \int_{Q_R} |D^m w|^p\, dx \right)^{1+
     \delta}.
\end{equation}

In the case $p \ge 2$ there holds for any $\mu \in [0,1]$:
\begin{equation}\label{HI.Dm+1.3}
     \int_{Q_{\tau R}} \left| D\left[ (\mu^2 + |D^m w|^2)^{\frac{p-2}{4}}D^m w\right]
     \right|^{2(1+\delta)}\, dx \le \left(\frac{c}{R^2} \int_{Q_R} (\mu^2 + 
     |D^m w|^2)^{p/2}\, dx\right)^{1+\delta}. 
\end{equation}
\hfill $\blacksquare$
\end{Lem}

\begin{proof}
Since such a result is more or less standard, we only show the main ideas of the proof here. Our aim is to show a reverse H\"older inequality which translates via
Gehring's lemma into the desired higher integrability result. Distinguishing
both the cases $p \ge 2$ and $p < 2$ and $\mu \not= 0$, $\mu = 0$, we proceed as 
follows: As in the proof of Lemma \ref{A priori Lemma}, we test system \eqref{Gl.Konst.} with the function $\phi \equiv D(\eta^{2m}
D(w-P))$ with a suitable cut off function $\eta$ and a polynomial of degree $m$ which
we specify later. 

In the case $p \ge 2$, following the estimates in the proof of Lemma \ref{A priori Lemma}, see \eqref{Absch.1}, finally applying Young's inequality,
we obtain
\begin{equation*}
\begin{aligned}
     \int_{Q_{\tau R}} &(\mu^2 + |D^m w|^2)^{\frac{p-2}{2}} |D^{m+1}w|^2\, dx\\
     &\le
     c \left[ \int_{Q_{\sqrt{\tau}R}} (\mu^2 + |D^m w|^2)^{p/2}\, dx +
     \sum_{k=1}^m \int_{Q_{\sqrt{\tau}R}} R^{-kp} |D^{m+1-k}(w-P)|^p\, dx \right].
\end{aligned}
\end{equation*}
Having in mind \eqref{CZ.Apr.est.techn.1}, the left hand side is estimated from below
by
\begin{equation*}
     \int_{Q_{\tau R}} \left| D\left[ (\mu^2 + |D^m w|^2)^{\frac{p-2}{4}}D^m w \right]
     \right|^2\, dx.
\end{equation*}
Choosing the coefficients of the polynomial $P$ in such a way that
\begin{equation*}
     \int_{Q_{\sqrt{\tau}R}} D^j(w-P)\, dx = 0, \qquad \mbox{ for }\ j = 0,\ldots,m-1
     \, ,
\end{equation*}
we can use Poincar\'e's inequality and elementary algebraic calculations to estimate
the second term of the right hand side from above by
\begin{equation*}
     cR^{-2} \int_{Q_{\sqrt{\tau}R}} (\mu^2 + |D^m(w-P)|^2)^{\frac{p-2}{2}}|D^m(w-P)
     |^2\, dx.
\end{equation*}
Now choosing the highest order coefficients of the polynomial $P$ such that
\begin{equation*}
\begin{aligned}
     \int_{Q_{\sqrt{\tau}R}} &(\mu^2 + |D^m(w-P)|^2)^{\frac{p-2}{2}}|D^m (w-P)|^2\, dx\\
     &= \int_{Q_{\sqrt{\tau}R}} \left| (\mu^2 + |D^m w|^2)^{\frac{p-2}{4}}D^m w - 
     \left( (\mu^2 + |D^m w|^2)^{\frac{p-2}{4}} D^m w\right)_{Q_{\sqrt{\tau}R}} 
     \right|^2\, dx,
\end{aligned}
\end{equation*}
we can apply Sobolev-Poincar\'e's inequality to conclude a reverse H\"older inequality
of the type
\begin{eqnarray*}
     &&\mskip-20mu \int_{Q_{\tau R}} \left| D\left[ (\mu^2 + |D^m w|^2)^{\frac{p-2}{4}}D^m w\right]
     \right|^2\, dx \\
     &&\mskip+15mu \le c \Biggl[\left( \int_{Q_{\sqrt{\tau} R}} \left| D\left[ (\mu^2
     + |D^m w|^2)^{\frac{p-2}{4}})D^m w \right]\right|^{\frac{2n}{n+2}}\, dx 
     \right)^{\frac{n+2}{n}}
     + R^{-2} \int_{Q_{\sqrt{\tau}R}} (\mu^2 + |D^m w|^2)^{p/2}\, dx \Biggr].
\end{eqnarray*}
Gehring's Lemma now provides the desired higher integrability. Combining this result
with the estimates in Lemma \eqref{A priori Lemma}, we end up with \eqref{HI.Dm+1.3}.
\begin{remark}
Here we also need higher integrability of $|D^m w|^{p}$, which is standard to prove. See
for example \cite{Giusti:2003} for higher integrability results of this type.
\end{remark}

In the case $1 < p < 2$ and $\mu \not= 0$, we finally obtain, using the same test function
as above (see \eqref{Absch. 4}) with a polynomial $P$ of degree $m-1$ and taking into consideration \eqref{CZ.Apr.est.diff}
\begin{eqnarray*}
     &&\int_{Q_{\tau R}}\left| D\left[(\mu^2 + |D^m w|^2)^{p/4}\right]\right|^2\, dx  \\
     && \mskip+130mu \le \frac{c}{R^2}\Biggl[\int_{Q_{\sqrt{\tau}R}} \left| (\mu^2 + |D^m w|^2)^{p/4} - \left( (\mu^2 + |D^m w|^2)^{p/4}\right)_{Q_{\sqrt{\tau}R}}\right|^2\, dx\\
     &&\mskip+230mu  + \int_{Q_{\sqrt{\tau}R}} (\mu^2 + |D^m w|^2)^{p/2}\, dx \Biggr].
\end{eqnarray*}
We now apply Sobolev-Poincar\'e's inequality, obtaining a reverse H\"older inequality of 
which allows us to apply Sobolev-Poincar\'e's inequality, obtaining a reverse H\"older
inequality of the type 
\begin{eqnarray*}
     &&\int_{Q_{\tau R}}\left| D\left[(\mu^2 + |D^m w|^2)^{p/4}\right]\right|^2\, dx\\
     &&\mskip+60mu \le c\left[\left( \int_{Q_{\sqrt{\tau}R}} \left|D\left[(\mu^2 + |D^m w|^2)^{p/4}\right]\right|^{\frac{2n}{n+2}}\, dx \right)^{\frac{n+2}{n}} + R^{-2}\int_{Q_{\sqrt{\tau}R}} (\mu^2 + |D^m w|^2)^{p/2}\, dx \right].
\end{eqnarray*}
Again applying Gehrings Lemma and combining the result with the estimate of Lemma
\ref{A priori Lemma} and a priori higher integrability for $|D^m w|^p$, provides the desired estimate \eqref{HI.Dm+1.1}.
In the case $\mu = 0$ we have by \eqref{Absch. 4} and \eqref{CZ.Apr.est.techn.1}
\begin{eqnarray*}
     &&\int_{Q_{\tau R}} \left| D\left[ |D^m w|^{\frac{p-2}{2}}D^m w\right]\right|^2\, dx \\
     &&\mskip+50mu \le \frac{c}{R^2} \Biggl[\int_{Q_{\sqrt{\tau}R}} \left||D^m w|^{\frac{p-2}{2}}D^m w - \left(|D^m w|^{\frac{p-2}{2}}D^m w\right)_{Q_{\sqrt{\tau}R}}\right|^2\, dx
     + \int_{Q_{\sqrt{\tau}R}} |D^m w|^p\, dx \Biggr].
\end{eqnarray*}
Sobolev-Poincar\'e's inequality now provides
\begin{eqnarray*}
     &&\int_{Q_{\tau R}} \left| D\left[ |D^m w|^{\frac{p-2}{2}}D^m w\right]\right|^2\, dx \\
     &&\mskip+100mu \le c \left[\left( \int_{Q_{\tau R}} \left| D\left[ |D^m w|^{\frac{p-2}{2}}D^m 
     w\right]\right|^{\frac{2n}{n+2}}\, dx\right)^{\frac{n+2}{n}} + R^{-2}\int_{Q_{\sqrt{\tau}R}} |D^m w|^p\, dx\right].
\end{eqnarray*}
Again Gehrings Lemma and finally the estimates of Lemma \ref{A priori Lemma} lead to the
desired estimate \eqref{HI.Dm+1.2}.
\end{proof}

\section{Calder\'on--Zygmund coverings}

We consider a cube $Q_0 \subset \bR^n$ and define by ${\mathcal D}(Q_0)$ the
set of all dyadic subcubes $Q$ of $Q_0$, i.e. those cubes with sides parallel to 
the sides of $Q_0$ that can be obtained from $Q_0$ by a positive finite number of 
dyadic subdivisions. We call $Q_p$ {\em  a predecessor} of $Q$, if $Q$ is obtained from $Q_p$ by a finite number of dyadic subdivisions. In particular we call 
$\tilde{Q} \in {\mathcal D}(Q_0)$ {\em the predecessor} of $Q$, if $Q$ is obtained 
from $\tilde{Q}$ by exactly one dyadic subdivision from $\tilde{Q}$.

The following lemma will play an essential role in the proof Theorem \ref{Hauptsatz 1}. The proof is done by Calder\'on--Zygmund coverings and can be found
for example in \cite{Caffarelli:1998}.

\begin{Lem}\label{Calderon-Zygmund Lemma}
Let $Q_0 \subset \bR^n$ be a cube. Moreover let $X\subset Y \subset Q_0$ be measurable 
sets satisfying the following: There exists $\delta > 0$ such that
\begin{itemize}
\item[(i)]
\begin{equation*}
     |X| < \delta |Q_0|,
\end{equation*}
\end{itemize}
and
\begin{itemize}
\item[(ii)] 
for any cube $Q \in {\mathcal D}(Q_0)$ there holds
\begin{equation*}
     |X \cap Q| > \delta |Q| \Longrightarrow \tilde{Q} \subset Y,
\end{equation*}
in which $\tilde{Q}$ denotes the predecessor of $Q$. 
\end{itemize}
Then there holds
\begin{equation*}
     |X| < \delta |Y|.
\end{equation*}
\strut\hfill$\blacksquare$
\end{Lem}

\section{Hardy Littlewood maximal function}

We will use properties of the Hardy Littlewood maximal function, which we will state here without proving them. For a more detailed discussion about maximal operators 
see \cite{Iwaniec:1992} and \cite{Stein:1993}. 

Let $Q_0 \subset \bR^n$ be a cube. For a function $f \in L^1\left(Q_0\right)$ we define the
restricted maximal function on $Q_0$ by 
\begin{equation}
     M^*_{Q_0}(f)(x) := \sup\limits_{Q \subseteq Q_0, x \in Q} \midint_Q 
     \left|f(y)\right|\, dy, \quad x \in Q_0,
\end{equation}
where $Q$ denotes an arbitrary subcube of $Q_0$, not necessarily centered 
in $x \in Q_0$. In an analogue way we define for $s > 1$ and $f \in
L^s\left(Q_0\right)$ 
\begin{equation}
     M^*_{s,Q_0}(f)(x) := \sup\limits_{Q \subseteq Q_0, x \in Q} \left(\midint_Q 
     \left|f(y)\right|^s\, dy\right)^{1/s}.
\end{equation}

We will need the following properties of the maximal function operator:
\begin{Lem}\label{Absch.Maxfkt.}
For $Q_0 \subset \bR^n$ and $s > 1$ let the maximal functions $M^*_{Q_0}$ and 
$M^*_{s,Q_0}$ be defined as above. Then the following estimates hold:
\begin{itemize}
\item[(M1)] For $f \in L^1\left(Q_0\right)$ and for any $\alpha > 0$ there holds
\begin{equation}\label{Absch.Mfkt.1}
     \left| \left\{ x: M^*_{Q_0}(f)(x) > \alpha \right\} \right| \le 
     \frac{c_W}{\alpha} \int_{Q_0} \left|f(y)\right|\, dy,
\end{equation}
with a constant $c_W \equiv c_W(n)$, for example $c_W = 3^n$ sufficies.
\item[(M2)] For $f \in L^p\left(Q_0\right), 1 < p < \infty$ we have $M^*(f) \in 
L^p$ and there holds
\begin{equation}\label{Absch.Mfkt.2}
     \int_{Q_0} \left|M^*_{Q_0}(f)(y)\right|^p\, dy \le \frac{3^ne p^2}{p-1} 
     \int_{Q_0} \left|f(y)\right|^p \, dy.
\end{equation}
\item[(M3)] An inequality similar to the one in (M2) holds also for the maximal
function $M^*_{s,Q_0}$, i.e. for $p > s$ we have:
\begin{equation}\label{Absch.Mfkt.3}
     \int_{Q_0} \left|M^*_{s,Q_0}(f)(y)\right|^p\, dy \le \frac{3^ne p^2}{s(p-s)} 
     \int_{Q_0} \left|f(y)\right|^p \, dy.
\end{equation}
\end{itemize}
\strut\hfill$\blacksquare$
\end{Lem}
A direct consequence of (M1) is the following
\begin{corollary}\label{Kor.Maxfkt.}
Let $f \in L^s(Q_0)$ with $s > 1$. Then 
\begin{equation*}
     \left| \left\{ x: M^*_{Q_0}(f)(x) > \alpha \right\} \right| \le 
     \frac{c_W}{\alpha^s} \int_{Q_0} \left|f(y)\right|^s\, dy,
\end{equation*}
with the constant $c_W$ of $(M1)$.
\strut\hfill$\blacksquare$
\end{corollary}

\section{Some technicalities}

The following technical lemmas will be used at several points in the proof of the main 
theorem. Since they are more or less standard we will only cite them. 

\begin{Lem}[\cite{Cupini:1999}, Lemma 2.2]\label{techn. 3}
Let $p > 1$. Then there exists a constant $c$, such that for every $\mu \ge 0$,
$\xi, \eta \in \bR^k$ there holds
\begin{equation*}
     (\mu^2 + |\xi|^2)^{p/2} \le c(\mu^2 + |\eta|^2)^{p/2} + c(\mu^2 + |\xi|^2 +
     |\eta|^2)^{(p-2)/2}|\xi - \eta|^2.
\end{equation*}
\strut \hfill $\blacksquare$
\end{Lem}

\begin{Lem}[see \cite{Campanato:1982}]\label{techn. Lemma}
Let $a,b \in \bR^N$ and $\nu > -1$. Then there exist constants $c(\nu)$, 
$C(\nu) > 0$, such that there holds
\begin{equation*}
     c(\nu)\left(\mu+ |a| + |b|\right)^{\nu} \le \int_0^1 \left(\mu+|a+tb|\right)^{
     \nu}\, dt \le C(\nu) \left(\mu+|a|+|b|\right)^{\nu}.
\end{equation*}
\strut\hfill $\blacksquare$
\end{Lem}

\section{Proof of the main theorem}

The proof of Theorem \ref{Hauptsatz 1} is at many stages similar or identical to the
proof in the case $m=1$, which is done in \cite{Acerbi:2004}. Therefore
some of the estimates will only be cited (for example the comparison estimate). We
will especially point out the differences to the higher order case here.

\subsection{Choice of constants and radii (I)}

To proceed with the proof of Theorem \ref{Hauptsatz 1} we initially fix some of the constants.
The proof will take place on the cube $Q_{4R_0} \Subset \Omega$.
The radius $R_0$ of this cube will be restricted at several points in the course of the
proof. At first we choose the radius so small that
\begin{equation}\label{Lokalisierung}
     \left\{ \begin{array}{ll} \omega\left(8nR_0\right) \le \sqrt{\frac{n+1}{
     n}} -1, \\[0.3cm] 0 < \omega(R) \log \left(\frac{1}{R}\right) \le \tilde{L}, \quad
     \mbox{ for any } 0 < R \le 8nR_0.\\ \end{array} \right.
\end{equation}
Furthermore we set
\begin{equation}\label{Def. K_0}
     K_0 := \int_{Q_{4R_0}} \left| D^m u\right|^{p(x)}\, dx + 1.
\end{equation}

\subsection{Higher integrability}

We will show that the condition
\begin{equation}\label{Stet. Mod. schw.}
     \lim\limits_{\rho \downarrow 0} \omega(R) \log\left(\tfrac{1}{\rho}\right) \le M 
     < + \infty,
\end{equation}
on the modulus of continuity $\omega$ yields a certain higher integrability for
$|D^m u|^{p(\cdot)}$. We note that condition \eqref{Stet. Mod. schw.} is weaker than condition (\ref{Stet. Mod.}) which is needed for proving the main theorem.
Our result is the following
\begin{Lem}[Higher integrability of $\left|D^m u\right|^{p(\cdot)}$]\label{High.Int.}
Let $u \in W^{m,1}\left(\Omega;\bR^N\right)$ with $\left|D^m u\right|^{p(\cdot)} \in 
L^1_{loc}\left(\Omega\right)$ be a weak solution of (\ref{Gleichung}) under the
conditions (\ref{Beschr. p}), (\ref{Stet. Mod. schw.}) as well as 
\begin{equation}\label{Schranke A}
     \left|A\left(x,z\right)\right| \le L \left(1+|z|^2\right)^{\left(p(x)-1
     \right)/2},
\end{equation}
and
\begin{equation}\label{Ellipt. A}
     \nu\left(\mu^2 + |z|^2\right)^{p(x)/2} - L \le \left< A\left(x,z\right),z
     \right>,
\end{equation}
for all $x \in \Omega, z \in \bR^{\cN}$. Moreover let $F \in L^{p(\cdot)q}(\Omega;\bR^{\cN})$ for some $q > 1$.
Then there exist constants $c \equiv c(n,N,m,\gamma_1,\gamma_2,L/\nu,\tilde{L}, 
M)$ and $c_g \equiv c_g(n,N,m,\gamma_1,\gamma_2,\nu,
L,\tilde{L})$ such 
that the following holds:
If $R_0$ is the radius from (\ref{Lokalisierung}), $K_0$ from \eqref{Def. K_0}, $Q_{4R_0} \Subset \Omega$, $\sigma > 0$ a constant with 
\begin{eqnarray*}
     0 < \sigma \le \sigma_0 := \min\left\{ \frac{c_g}{K_0^{\frac{2q\omega\left( 
     8nR_0\right)}{\gamma_1}}},\ q-1,\ 1 \right\},
\end{eqnarray*}
then for every cube $Q_R \subseteq Q_{4R_0}$ we have that
\begin{equation}\label{High.Int.Res.}
     \Biggl( \midint_{Q_{R/2}} \mskip-20mu \left|D^m u\right|^{p(x)\left(1 
          +\sigma\right)} \, dx\Biggr)^{\frac{1}{1+\sigma}}
     \le c \midint_{Q_R} \mskip-15mu \left|D^m u \right|^{p(x)}\, dx +
          c \Biggl( \midint_{Q_R} \mskip-15mu \left|F\right|^{
     p(x)\left(1+\sigma\right)}\, dx + 1 \Biggr)^{\frac{1}{1+\sigma}}. 
\end{equation}
\strut\hfill $\blacksquare$
\end{Lem}

\begin{proof}
Since the proof of this result is in many points similar to the proof in 
the case $m=1$, we only show the main steps here. 
Let $Q_R \subseteq Q_{4R_0}$ be a cube and 
\begin{equation}
     p_1 := \inf\left\{p(x): x \in Q_R\right\},\ \ \ \ \ p_2 := \sup
     \left\{ p(x): x \in Q_R \right\}.
\end{equation}
Then, $p_2 - p_1 \le \omega\left(2nR\right)$ and by the choice of $R_0$ in  (\ref{Lokalisierung}) we have
\begin{equation}
\frac{p_2}{p_1} = \frac{p_2-p_1}{p_1} + 1 \le \sqrt{\frac{n+1}{n}} =: \s.
\end{equation}
We test system \eqref{Gleichung} by the function $\phi_1 \equiv \eta^{mp_2} \left(u-P\right)$, where $\eta \in C^{\infty}_c\left(Q_R\right)$ 
denotes a standard cut-off function satisfying $0 \le \eta \le 1$ and $\eta \equiv 1$ 
on $Q_{R/2}$ as well as $\left|D^k \eta\right| \le \frac{1}{R^k}$ for $k=1,\ldots,m$ 
and $P:\bR^n \to \bR$ denotes the unique polynomial of degree $m-1$ satisfying
\begin{equation}\label{Polynom}
     \left( D^k \left(u - P\right) \right)_{Q_R} = \midint_{Q_R} D^k \left(u-
     P\right) \,dx = 0 \ \ \mbox{ for }\ \ k=0,\ldots,m-1. 
\end{equation}
It is easy to see that we have
\begin{equation}\label{Abl.Eta}
     \left|D^k\left(\eta^{mp_2}\right)\right| \le C\left(n,m,\gamma_2\right) 
          R^{-k}\sum\limits_{j=1}^k \eta^{mp_2-j}
     \le C\left(n,m,k,\gamma_2\right) R^{-k} \eta^{mp_2 - k}. 
\end{equation}
Setting in the test function and using \eqref{Ellipt. A} we obtain
\begin{eqnarray*}
     \nu\int_{Q_R} \eta^{mp_2} \left|D^m u\right|^{p(x)}\, dx - L
     &\le&\int_{Q_R} \left<A\left(x,D^m u\right),D^m \phi_1\right>\, dx\\
     &&\mskip+10mu - \int_{Q_R} \sum\limits_{k=1}^m \binom{m}{k}\left<A\left(
          x,D^m u\right), D^k\left(\eta^{mp_2}\right) \odot D^{m-k} 
          \left(u-P\right) \right>\, dx\\ 
     &=& I_1 + I_2,
\end{eqnarray*}
with the obvious labelling. 
Using \eqref{Gleichung}, $I_1$ can be estimated by applying Young's inequality several times in a standard way (note that the constant in Young's
inequality may depend on $p(x)$; writing the constant down explicitely, one can
easily see that it can be estimated by a constant depending only on $\gamma_1$ 
and $\gamma_2$) and using the fact that $p(x) \le p_2$ on $Q_R$ to obtain
\begin{equation}
\begin{aligned}
     I_1 \le &\ \varepsilon \int_{Q_R} \eta^{mp_2} \left|D^m u\right|^{p(x)}
          \, dx\\ 
     &+ c \left[ \int_{Q_R} \left|F(x)\right|^{p(x)}\, dx
      + \sum\limits_{k=1}^m 
          \int_{Q_R}\left(\frac{\left|D^{m-k}\left(u-P\right)\right|^{p_2}}{R^{kp_2
          }}+1 \right)\, dx \right],
\end{aligned}
\end{equation}
where $c \equiv c\left(m,\varepsilon,\gamma_1,\gamma_2\right)$.
By \eqref{Schranke A}, \eqref{Abl.Eta} and Young's inequality we estimate
\begin{equation*}
     I_2 \le  cL\sum\limits_{k=1}^m
          \int_{Q_R}\mskip-10mu\left(1+\left|D^m u\right|^{p(x)-1}\right) 
          \eta^{mp_2-k} \frac{|D^{m-k}(u-P)|}{R^k}\, dx,
\end{equation*}
where $c \equiv c(n,m,\gamma_2)$. 
We rewrite the integral appearing on the right hand side as follows :
\begin{eqnarray*}
     &&\int_{Q_R} \left(1+|D^m u|^{p(x)-1}\right)\eta^{mp_2-k} \frac{|D^{m-k}(u-P)|}{
          R^k}\, dx\\
     &&\mskip+60mu = \int_{Q_R} |D^m u|^{p(x)-1}\eta^{mp_2-k}\frac{|D^{m-k}( 
          u-P)|}{R^k}\, dx
          + \int_{Q_R} \eta^{mp_2-k}\frac{|D^{m-k}(u-P)|}{R^k}\, dx\\
     &&\mskip+60mu = I_{2,1}^{(k)} + I_{2,2}^{(k)}.
\end{eqnarray*}
The second integral we treat as usual noting that 
\begin{equation*}
     I_{2,2}^{(k)} \le \int_{Q_R} \frac{|D^{m-k} (u-P)|^{p_2}}{R^{kp_2}} + 1\, dx.
\end{equation*}
Therefore it remains to get a bound for $I_{2,1}^{(k)}$.
Since $\eta \le 1$ we have $\eta^{mp_2-k} \le 
\eta^{m(p_2-1)}$ for $k=1,\ldots,m$, and hence by Young's inequality and
$\frac{p(x)}{p(x)-1} \ge \frac{p_2}{p_2-1}$ we obtain that
\begin{equation*}
     I_{2,1}^{(k)} 
     \le  \varepsilon \int_{Q_R} \eta^{mp_2} |D^m u|^{p(x)}\, dx + 
          c \int_{Q_R} \frac{|D^{m-k}(u-P)|^{p_2}}{R^{kp_2}} + 1\, dx.
\end{equation*}
Combining the estimates for $I_{2,1}^{(k)}$ and $I_{2,2}^{(k)}$ we finally arrive at
\begin{eqnarray*}
     \nu \int_{Q_R} \eta^{mp_2}|D^m u|^{p(x)}\, dx
     &\le& cL\varepsilon \int_{Q_R} \eta^{mp_2}|D^m u|^{p(x)}\, dx\\
     &&+ cL \Biggl[ \int_{Q_R}|F(x)|^{p(x)}\, dx
          + \sum\limits_{k=1}^m \int_{Q_R} \left(\frac{
           |D^{m-k}(u-P)|^{p_2}}{R^{kp_2}}+1\right)\, dx  \Biggr],
\end{eqnarray*}
where the constant $c$ depends only on $n,m,\gamma_1,\gamma_2$ and $\varepsilon$.
Now choosing as usual $\varepsilon = \frac{\nu}{2cL}$ (note that we can also assume
that $\varepsilon < 1 - 1/\gamma_2$ by choosing $c$ large enough) we can absorb
the first integral on the right hand side. 
Dividing the resulting inequality by $\nu/2$ leads us 
to
\begin{equation}\label{High.Int.Cacc.}
     \int_{Q_{R/2}} |D^m u|^{p(x)}\, dx 
     \le c \int_{Q_R} |F(x)|^{p(x)}\, dx + c \sum\limits_{k=1}^m \int_{Q_R} 
     \left(\frac{|D^{m-k}(u-P)|^{p_2}}{R^{kp_2}} + 1\right)\, dx,
\end{equation}
where $c$ depends only on $n,m,\gamma_1,\gamma_2$ and $L/\nu$.
Taking into account the properties (\ref{Polynom}) of the polynomial $P$, using Poincar\'e's inequality and taking the mean values on both sides, we arrive at:
\begin{equation*}
     \midint_{Q_{R/2}} \left|D^m u\right|^{p(x)}\, dx \le c \midint_{Q_R}\left| 
     F(x)\right|^{p(x)}\, dx + c \midint_{Q_R} \left( \left|\frac{D^{m-1} 
     \left(u-P\right)}{R}\right|^{p_2}+1\right)\, dx.
\end{equation*}
Taking into account the definition of $P$, i.e. (\ref{Polynom}), we can apply Sobolev-Poincar\'e's 
inequality on the right hand side with exponents $p \equiv p_2$ and 
$p_* \equiv \frac{np_2}{n+p_2}$ to obtain 
\begin{equation*}
     \midint_{Q_R} \left|\frac{D^{m-1}(u-P)}{R}\right|^{p_2}\, dx 
     \le c_{SP} \left( \midint_{Q_R}|D^m u|^{\frac{np_2}{n+p_2}}\, dx \right)^{ 
     \frac{n+p_2}{n}}.
\end{equation*}
Note that the constants
in Poincar\'e's and Sobolev-Poincar\'e's inequalities can be replaced by constants that only depend on $\gamma_2$ instead of $p_2$, thus $c_{SP} \equiv c(n,N,\gamma_2)$.
H\"older's inequality, applied with the exponents $p \equiv \frac{p_1 
s(n+p_2)}{p_2(n+1)}> 1$ provides (note $\s^2 = \frac{n+1}{n}$)
\begin{equation*}
     \midint_{Q_R} \mskip-8mu |D^m u|^{\frac{np_2}{n+p_2}}\, dx \le 
     \left(\midint_{Q_R} \mskip-8mu |D^m u|^{p_1/\s}\, 
     dx\right)^{\frac{p_2(n+1)}{p_1\s(n+p_2)}}.
\end{equation*}
Therefore we obtain
\begin{equation*}
     \midint_{Q_{R/2}} |D^m u|^{p(x)}\, dx 
     \le c\midint_{Q_R} |F(x)|^{p(x)}\, dx +
           c \left(\midint_{Q_R} |D^m u|^{p(x)/\s}\, dx +1\right)^{\frac{p_2\s}{p_1}}.
\end{equation*}
Noting that $\frac{p_2 \s}{p_1} = \s + \s\left(\frac{p_2}{p_1}-1\right)$ and
$\frac{p_2}{p_1} -1 \le \frac{\omega(2nR)}{p_1}$ we get
\begin{equation*}
     \left(\midint_{Q_R} |D^m u|^{p(x)/\s}+1\, dx\right)^{\frac{p_2\s}{p_1}}
     \le c\left[1+\left(\midint_{Q_R} |D^m u|^{p(x)/\s}\, dx\right)^{ 
          \frac{\s \omega(2nR)}{p_1}} \left(\midint_{Q_R} |D^m u|^{p(x)/\s}\, dx 
          \right)^{\s}\right],
\end{equation*}
where $c \equiv c(n,\gamma_1,\gamma_2)$.
With $\frac{\s}{p_1} \le 2$ and $R \le 1$ we obtain for the first term on the right 
hand side
\begin{equation*}
     \left( \midint_{Q_R} |D^m u|^{p(x)/\s}\, dx\right)^{\frac{\s \omega(
     2nR)}{p_1}} 
     \le R^{-2n\omega(2nR)} \left(\int_{Q_R} |D^m u|^{p(x)}\, dx + |Q_R|
     \right)^{\frac{\s\omega(2nR)}{p_1}}.
\end{equation*}
Noting that by the localization properties \eqref{Lokalisierung} we see that
$R^{-2n\omega(2nR)} \le c(n,\tilde{L})$ and taking into account the definition of
$K_0$ we conclude the following reverse H\"older inequality (note also that
$K_0 > 1$ and $p_1 \ge \gamma_1$):
\begin{equation}\label{Rev. Hoelder}
     \midint_{Q_{R/2}} |D^m u|^{p(x)}\, dx
     \le cK_0^{\frac{2\omega(8nR_0)}{\gamma_1}} \left(\midint_{Q_R} |D^m u|^{p(x)/\s}
     \, dx\right)^{\s} + c\midint_{Q_R} \left(|F(x)|^{p(x)}+1\right)\, dx,
\end{equation}
where $K_0 > 1$ is from \eqref{Def. K_0} and $c \equiv c\left(n,N,m,M,\gamma_1, \gamma_2,L/\nu 
\right)$.
This inequality holds for any cube $Q_R \subseteq Q_{4R_0}$ and the appearing 
constants do not depend on the choice of the particular cube $Q_R$. Gehring's Lemma 
in the version which is written in \cite{Acerbi:2001b} with $f  \equiv \left|D^m 
u\right|^{p(x)/\s}$ and $\phi \equiv (\left|F\right|^{p(x)} + 1)^{1/\s}$ under 
consideration of the restriction on $\sigma$ finally provides the assertion.
\end{proof}

\subsection{Choice of constants and radii (II)}

First we observe that, since $K_0 \ge 1$ (see the definition of $K_0$ in 
\eqref{Def. K_0}), we have for any $K \ge K_0$:
\begin{equation}\label{Absch.sigma_0}
     \sigma_0 \ge \min\left\{1,q-1,c_g\right\} K^{-\frac{2q\omega\left(8nR_0\right)}{
     \gamma_1}},
\end{equation}
where $\sigma_0$ is the constant from Lemma \ref{High.Int.}.
We set
\begin{equation}\label{Def. K_M}
     K_M := \int_{\Omega} \left(\left|D^m u\right|^{p(x)} + |F|^{p(x)q} + 2\right)
     \, dx + 1,
\end{equation}
and
\begin{equation}\label{Def. sigma_M}
     \sigma_m := \min \left\{ \frac{c_g}{K_M^{\frac{2q\left(\gamma_2-\gamma_1\right)}{
     \gamma_1}}} ,\frac{q-1}{2},1\right\} > 0,\qquad \sigma_M := c_g + q.
\end{equation}
Therefore $K_M \ge K_0$.
Furthermore for any $1 \le K \le K_M$ we have 
\begin{equation}
     \sigma_m \le \sigma_0 \le \sigma_M.
\end{equation}
We now choose the higher integrability exponent $\sigma$ in Lemma \ref{High.Int.}
such that
\begin{equation}\label{Einschr.sigma}
     \sigma := \tilde{\sigma}\sigma_0 \ \mbox{ with }\ \ 0 < \tilde{\sigma} < \min
     \left\{ \gamma_1-1, 1/2\right\}.
\end{equation}
Then by (\ref{Absch.sigma_0}) we have for any $\beta \in \left[\frac{\gamma_2}{
\gamma_2-1}, \frac{\gamma_1}{\gamma_1-1}\right]$ and $K \ge K_0$:
\begin{equation}\label{Absch.sigma_2}
     \sigma^{-\beta} \le c \tilde{\sigma}^{-\beta} K^{\beta\frac{2q\omega\left(8nR_0
     \right)}{\gamma_1}} \le c\left(n,N,m,\gamma_1,\gamma_2,L/\nu,q\right) \tilde{
     \sigma}^{-\beta} K^{\frac{2q\omega\left(8nR_0\right)}{\gamma_1-1}}.
\end{equation}
By the choice of $\sigma$ in \eqref{Einschr.sigma} and the structure of the 
constant $\sigma_0$ in Lemma \ref{High.Int.} we have that 
\begin{equation}\label{Absch.sigma}
     \sigma < \frac{q-1}{2}.
\end{equation}
Now we impose for a fixed choice of $\tilde{\sigma}$ a further restriction on
the size of $R_0$ by claiming
\begin{equation}\label{Einschr.R_0}
     \max\left\{ 2q\omega\left(8nR_0\right),\frac{2q\omega\left(8nR_0\right)}{
     \gamma_1-1} \right\} \le \frac{\tilde{\sigma} \sigma_m}{4}.
\end{equation}
Therefore $R_0$ depends on $n,N,m,\gamma_1,\gamma_2,L/\nu,q, \| |D^m u(\cdot)|^{
p(\cdot)}\|_{L^1\left(\Omega\right)}, \| |F(\cdot)|^{p(\cdot)}\|_{L^1\left(\Omega\right)}$ and
$\tilde{\sigma}$. (\ref{Einschr.R_0}) immediately implies
\begin{equation}\label{Einschr.R_0.2}
     \omega\left(8nR_0\right) \le \max\left\{2q\omega\left(8nR_0\right), \frac{2q
     \omega\left(8nR_0\right)}{\gamma_1-1}\right\} \le \frac{\tilde{\sigma} \sigma_m}{
     4} \le \frac{\tilde{\sigma} \sigma_0}{4} = \frac{\sigma}{4}.
\end{equation}

\subsection{Calder\'on-Zygmund type estimate}\label{chap:CZ}

The key to the proof of Theorem \ref{Hauptsatz 1} is the following lemma, which is
an application of Lemma \ref{Calderon-Zygmund Lemma} to special sets $X$ and $Y$.

\begin{Lem}\label{Calderon-Zygmund}
Let $u \in W^{m,p(\cdot)}\left(\Omega;\bR^N\right)$ be a weak solution of system 
(\ref{Gleichung}) under the structure conditions 
(\ref{Wachst. A}), (\ref{Stetigk. A}),
(\ref{Beschr. p}) and (\ref{Stet. Mod.}). Furthermore let $\lambda \ge 1$, $0 < 
\tilde{\sigma} < 1$ as in (\ref{Einschr.sigma}) and $B_M > 1$. Then there exists a constant $A \equiv A(n,N,m,\gamma_1,\gamma_2,L/\nu,\tilde{L}) \ge 
2$, independent of $\lambda, \tilde{\sigma}, u, A, F, B_M$ and a radius 
\begin{equation*}
     R_1 \equiv R_1(n,N,m,\gamma_1,\gamma_2,L/\nu,\tilde{L},q,\tilde{\sigma},B_M),
\end{equation*}
such that the following holds:
If $R_0 \le R_1$ is so small that
(\ref{Lokalisierung}) and (\ref{Einschr.R_0}) hold, if $K_0,\sigma_0$ are the constants from
(\ref{Def. K_0}), resp. Lemma \ref{High.Int.}, if $\sigma := \tilde{\sigma} \sigma_0$ as in 
\eqref{Einschr.sigma}, if
\begin{equation}\label{Def. K}
     K := \int_{Q_{4R_0}} \left|D^m u\right|^{p(x)} + \left|F\right|^{p(x)(1+\sigma)}
     \, dx + 1,
\end{equation}
if $K_M$ and $\sigma_M$ are the constants from (\ref{Def. K_M}), resp. (\ref{Def. sigma_M}), then for
any $1 < B < B_M$ and $\vartheta \ge \bigl( B^{\frac{n(1+\delta)}{n-2}} \bigr)^{-1}$ there exists 
$\tilde{\varepsilon} > 0$, independent of $\lambda$, such that the following implication holds:

If on some $Q \in {\mathcal D}\left(Q_{R_0}\right)$ we have
\begin{equation}\label{Calderon-Zygmund 1}
\begin{aligned}
     \Bigl| Q \cap \Bigl\{ x\in Q_{R_0} :& M^*_{Q_{4R_0}} \left(\left|D^m 
          u(\cdot)\right|^{p(\cdot)}\right)(x) > ABK^{\sigma}\lambda,\\
     & M^*_{1+\sigma,Q_{4R_0}}\left(\left|F(\cdot)\right|^{p(\cdot)} + 1\right)( 
          x) < \tilde{\varepsilon}\lambda\quad \Bigr\} \Bigr| > \vartheta \left|Q 
          \right|,
\end{aligned}
\end{equation}
then for the predecessor $\tilde{Q}$ of $Q$ there holds
\begin{equation}\label{Calderon-Zygmund 2}
     \tilde{Q} \subseteq \left\{ x \in Q_{R_0}\ :\ M^*_{Q_{4R_0}} \left(\left|D^m u
     (\cdot)\right|^{p(\cdot)}\right)(x) > \lambda\right\}.
\end{equation}
\strut\hfill $\blacksquare$
\end{Lem}

\begin{proof}
We will prove the statement by contradiction. The constants $A, \tilde{\varepsilon}$
as well as the radius $R_1$ will be chosen at the end of the proof.
Let us assume that (\ref{Calderon-Zygmund 1}) holds, but (\ref{Calderon-Zygmund 2}) 
is false. Then there exists a point $x_0 \in \tilde{Q}$, such that 
\begin{equation*}
     M^*_{Q_{4R_0}} \left(\left|D^m u(\cdot)\right|^{p(\cdot)} \right)(x_0) \le \lambda,
\end{equation*}
i.e. we have
\begin{equation}\label{CS 3}
     \midint_C \left|D^m u(x)\right|^{p(x)}\, dx \le \lambda,
\end{equation}
for all cubes $C \subseteq Q_{4R_0}$ with $x_0 \in C$. We define $S := 
2 \tilde{Q}$. Since the cube $\tilde{Q}$ 
is obtained from the cube $Q_{R_0}$ by at least one dyadic subdivision, we have
$\tilde{Q} \subseteq Q_{R_0}$ and therefore $S \subseteq Q_{2R_0}$. 
Therefore by the smallness condition (\ref{Lokalisierung}) imposed on the radius $R_0$ there holds
\begin{equation}
     s := \mbox{diam}(2S) \le 8nR_0,\ \mbox{ and therefore }\ \omega(s) \le \sigma/4.
\end{equation} 
In particular, since by $x_0 \in 2S$ the cube $2S \subseteq Q_{4R_0}$ is an admissible cube in
the maximal function $M^*_{Q_{4R_0}}$, by (\ref{CS 3}) there holds
\begin{equation}\label{CS 4}
     \midint_{2S} \left|D^m u(x)\right|^{p(x)}\, dx \le \lambda.
\end{equation}
Additionally (\ref{Calderon-Zygmund 1}) implies
\begin{equation}
     \left| \left\{ x \in Q:\ M^*_{1+\sigma,Q_{4R_0}} \left(|F(\cdot)|^{p(\cdot)}+1\right)(x) < 
     \tilde{\varepsilon} \lambda \right\}\right| > 0,
\end{equation}
so that there exists at least one point $x \in Q$, in which the maximal function  
$M^*_{1+\sigma,Q_{4R_0}}$ of $|F(\cdot)|^{p(\cdot)}+1$ is small. Since $Q \subset 2S \subset Q_{4R_0}$, 
this implies
\begin{equation}\label{CS 5}
     \left(\midint_S \left(|F|^{p(x)}+1\right)^{1+\sigma}\, dx\right)^{\frac{1}{
     1+\sigma}} < \tilde{\varepsilon}\lambda,\ \ \left(\midint_{2S} \left(|F|^{p(x)}+1
     \right)^{1+\sigma}\, dx\right)^{\frac{1}{1+\sigma}} < \tilde{\varepsilon}\lambda.
\end{equation}
We now use the localization argument from before in order to estimate $p(x)$ in a point $x$ by 
constant exponents $p_1,p_2$. For this purpose we let
\begin{equation}
     p_1 := p(x_m) = \min\limits_{\overline{2S}} p(x),\ \ \ p_2 := p(x_M) = 
     \max\limits_{\overline{2S}} p(x),\ \ \ x_M,x_m \in \overline{2S}.
\end{equation}
Obviously the exponents $p_1$ and $p_2$ depend on the local situation, especially on
the cube $Q \in {\mathcal D}(Q_R)$. Thus in the following estimates it will
be necessary to take care of the dependencies of the occurring constants on $p_1$ and
$p_2$, eventually replacing them by constants which only depend on the global bounds
$\gamma_1$ and $\gamma_2$ for $p$.
We first use the fact that $2S \subseteq Q_{4R_0}$, obtaining by the restriction (\ref{Einschr.R_0}) on the radius $R_0$ for any $x \in \overline{2S}$
\begin{equation}\label{p(x) 1}
\begin{aligned}
     p_2 & =  (p_2-p_1)+p_1
         \le \omega\left(|x_M - x_m|\right) + p_1
         \le \omega(s) + p_1
         \le p_1\left(1+\omega(s)\right)\\
         &\le p(x)\left(1+\omega(s)\right)
         \le p(x)\left(1+\omega(s)+\sigma/4\right)
         \le p(x)\left(1+\sigma\right).
\end{aligned}
\end{equation}
By (\ref{Einschr.sigma}) we have $\tilde{\sigma} < \gamma_1-1$. Recalling that $\sigma = 
\tilde{\sigma} \sigma_0$ and $\sigma_0 < 1$ (see Lemma \ref{High.Int.}) and $\gamma_1 \le p_1$ we have $\sigma \le p_1 - 1$, which implies
\begin{equation}\label{p(x) 2}
\begin{aligned}
     p_2\left(1+\sigma/4\right) 
          &\le \left(p_1 + \omega(s)\right)\left(1+\sigma/4\right) 
          =   p_1 + p_1 \sigma/4 + \omega(s)\left(1+\sigma/4\right) \\
          &\le p_1\left(1+\sigma/4 + \omega(s)\right) 
          \le p(x)\left(1+\sigma/4 + \omega(s)\right)
          \le p(x)\left(1+\sigma\right). 
\end{aligned}
\end{equation}

\subsubsection{Higher integrability}\label{ch:CZ.CZ1}

First we note that the higher integrability from Lemma \ref{High.Int.} together with
\eqref{CS 5} leads to an upper bound for the integral $\midint_S |D^m u|^{p_2}\, dx$.
Since to prove this, we can follow line by line the estimates in \cite[(59), p 132]{Acerbi:2004}, we do not rewrite the estimates here. We obtain
\begin{equation}
     u \in W^{m,p_2}\left(S\right),
\end{equation}
and 
\begin{equation}\label{High.Int.u}
     \midint_{S} \left|D^m u\right|^{p_2}\, dx \le cK^{\sigma/4}\lambda,\qquad
     \int_S \left|D^m u\right|^{p_2}\, dx \le cK^{1+\sigma/4},
\end{equation}
where the constants depend only on $n,N,m, \gamma_1, \gamma_2,M$ and $L/\nu$.

\subsubsection{The frozen system}

We consider the Dirichlet problem
\begin{equation}\label{Eingefroren}
     \left\{ \begin{array}{lcl}
               \displaystyle{\int_S\left< A\left(x_M,D^m w\right),D^m \phi\right>\, 
               dx = 0} && \mbox{for all } \ \phi \in W^{m,p_2}_0(S;\bR^N)\\[0.5cm]
               w \equiv u && \mbox{on} \ \partial S\\
             \end{array} 
     \right. .
\end{equation} 

Since the vector field $A$ is frozen in the point $x_M$, where the exponent $p(x_M) = p_2$ is
constant, the structure conditions for the original vector field $A(x,z)$ (see (\ref{Wachst. A}) and the remark after \eqref{Stet. Mod.}) lead to the following structure conditions for the frozen vector field $A(x_M,z)$:
\begin{eqnarray}
     &&\nu \left(\mu^2 + |z|^2\right)^{(p_2-2)/2}|\lambda|^2 \le \left< D_z A(x_M,z) \lambda
          , \lambda \right>
          \le L \left(\mu^2 + |z|^2\right)^{(p_2-2)/2} |\lambda|^2,\label{Str. gefr. 1a}\\[0.3cm]
     &&|A\left(x_M,z\right)|\le L\left(\mu^2+|z|^2\right)^{(p_2-1)/2},
          \label{Str. gefr. 2a}
\end{eqnarray}
for all $x \in S, z \in \bR^{\cN}$, where $\mu \in [0.1]$. It is easy to see that one can write these conditions also in the form
\begin{eqnarray}
      &&\mskip+30mu \nu \left(\mu^2 + |z_1|^2 + |z_2|^2\right)^{(p_2-2)/2} 
           |z_2-z_1|^2\le \left< A\left(x_M,z_2\right) - A\left(x_M,z_1\right),z_2-z_1 
           \right>\label{Str. gefr. 1},\\[0.3cm]
      &&\mskip+30mu \tfrac{\nu}{2} |z|^{p_2} \le \left< A\left(x_M,z\right),z\right> 
           + \nu \left(\left(L/ \nu\right)^{p_2}+1\right), \label{Str. gefr. 3}
 \end{eqnarray}
for all $z,z_1,z_2 \in \bR^{\cN}$.

The Dirichlet problem (\ref{Eingefroren}) admits a unique solution  
$w \in u + W^{m,p_2}_0(S;\bR^N)$.
 

Since the vector field $A(x_M,z)$ fulfills the hypothesis of Lemma \ref{A priori Lemma} (with $p,Q_{R}$
replaced by $p_2,S$), we can apply the lemma in combination with Lemma \ref{A priori Gehring} to conclude that 
\begin{itemize}
\item[$\bullet$] in the case $p_2 \ge 2$ we have the estimate
     \begin{equation}\label{Est.w.p_2>2}
          \int_{\frac{3}{4}S} \left| DV_{\mu}(D^m w) \right|^{2(1+\tilde{\delta})}
          \, dx \le
          \left(\frac{c}{R^2} \int_{S} \left|H_{\mu}(D^m w)\right|^2\, dx\right)^{ 
          1+\tilde{\delta}},
     \end{equation}
for any $\mu \in [0,1]$,
\item[$\bullet$] in the case $1 < p_2 < 2$, for $\mu \in (0,1]$ there holds
     \begin{equation}\label{Est.w.p_2<2.1}
          \int_{\frac{3}{4}S} \left| DH_{\mu}(D^m w) \right|^{2(1+\tilde{\delta})}
          \, dx \le
          \left( \frac{c}{R^2} \int_S \left|H_{\mu}(D^m w) \right|^2\, dx\right)^{1+\tilde{\delta}},
     \end{equation}
and for $\mu = 0$ we have
     \begin{equation}\label{Est.w.p_2<2.2}
          \int_{\frac{3}{4}S} \left| DV_0(D^m w)\right|^{2(1+\tilde{\delta})}
          \, dx \le \left( \frac{c}{R^2} \int_S \left|H_0(D^m w)\right|^2
          \, dx\right)^{1+\tilde{\delta}}.
     \end{equation}
\end{itemize}
with $\tilde{\delta} \equiv \tilde{\delta}(n,m,N,\gamma_1,\gamma_2,L,\nu) > 0$. Here we used the abbrevitations
\begin{equation}\label{Def.Vmu,Hmu}
     V_{\mu}(D^m w) := \left(\mu^2 + |D^m w|^2\right)^{\frac{p_2-2}{4}}D^m w\ , \qquad
     H_{\mu}(D^m w) := \left(\mu^2 + |D^m w|^2\right)^{p_2/4}.
\end{equation}
Note that the function $V_\mu:\bR^k \to \bR^k$ is quite common in the recent papers about regularity for systems and functionals. We will use here the following algebraic
property of $V_\mu$: For any $\mu \in [0,1], z,\eta \in \bR^k$ there exists a
constant $c \equiv c(n,\gamma_1,\gamma_2)$ such that 
\begin{equation}\label{prop.V}
     c^{-1}|z-\eta|(\mu^2 + |z|^2 + |\eta|^2)^{\frac{p_2-2}{4}} \le |V_{\mu}(z) -
     V_{\mu}(\eta)| \le c|z-\eta| (\mu^2 + |z|^2 + |\eta|^2)^{\frac{p_2-2}{4}}.
\end{equation}
The proof of this property can be found for instance in \cite{Acerbi:1989, Hamburger:1992}.
\subsubsection{Energy estimate}

Now we want to show the following energy estimate:
\begin{equation}\label{Energ.Absch.}
     \midint_S |D^m w|^{p_2}\, dx \le c\left(\gamma_2,L/\nu\right)
     \midint_S \left(|D^m u|^{p_2} + 1\right)\, dx.
\end{equation}
To prove \eqref{Energ.Absch.} we test (\ref{Eingefroren}) with $\phi = 
u-w$. $\phi$ is an admissible test function, since $u,w \in W^{m,p_2}$ and $w-u \in W^{m,p_2}_0\left(S;
\bR^N\right)$. Using (\ref{Str. gefr. 3}) we obtain
\begin{eqnarray*}
     \nu/2 \int_S |D^m w|^{p_2}\, dx 
     &\le& \int_S \left< A(x_M,D^m w),D^m w\right>\, dx + \nu\left( (L/\nu)^{p_2}+1
           \right)\\
     &=  & \int_S\left<A(x_M,D^m w),D^m u\right>\, dx +\nu\left( (L/\nu)^{p_2}+1
           \right)\\
     &\le& \int_S |A(x_M,D^m w)| |D^m u|\, dx +\nu\left( (L/\nu)^{p_2}+1
           \right) .
\end{eqnarray*}
The growth assumption (\ref{Str. gefr. 2a}) gives by Young's inequality
\begin{eqnarray*}
     \int_S |A(x_M,D^m w)| |D^m u|\, dx 
     &\le&L \int_S\left(1+|D^m w|^2\right)^{\frac{p_2-1}{2}} |D^m u|
           \, dx\\
     &\le&\varepsilon 2^{p_2/2-1} \int_S |D^m w|^{p_2}\, dx + 
           \varepsilon 2^{p_2/2-1} + \varepsilon^{1-p_2}L^{p_2} \int_S |D^m 
           u|^{p_2}\, dx.
\end{eqnarray*}
Combining these estimates, the asserted estimate follows by choosing
$\varepsilon = 2^{1-p_2/2}\nu/4$.

\subsubsection{Comparison estimate}

The next step is to establish a comparison estimate between $D^m u$ and $D^m w$. 
This turns out to be quite complicated, involving all the results from before, for
example the energy estimate, higher integrability, the structure conditions of the frozen system, the localization and very fine estimates on the $L \log^{\beta}L$ 
scale. Nevertheless the argument and estimates are nearly the same as in the
second order case and can be taken from \cite[p 134ff]{Acerbi:2004}. Note that at this point the continuity assumption \eqref{Stetigk. A} comes into play. Although this condition differs from the condition in \cite{Acerbi:2004}, the proof of the comparison estimate can be left unchanged. The reason for this is that we compare the original problem to a problem which is frozen in the point $x_M$ where the growth exponent is the maximal exponent $p(x_M) = p_2$. Therefore in the application of \eqref{Stetigk. A} we pass over from the exponents $p(x)$ and $p(y)$ on the right hand side to this maximal value $p_2$. Since this is the only point where the continuity condition comes into play, there is no problem in replacing the continuity condition in \cite{Acerbi:2004} by our assumption \eqref{Stetigk. A}.

Finally, one ends up with 
\begin{eqnarray}
     &&\midint_S (\mu^2 + |D^m u|^2 + |D^m w|^2)^{\frac{p_2-2}{2}}|D^m u - D^m w|^2
          \, dx \label{Vgl.Absch.}\\
     &&\mskip+100mu \le c_1\omega(s)\log\left(\frac{1}{s}\right) K^{\sigma} \lambda 
          + c_2\omega(s) \tilde{\sigma}^{-1} K^{\sigma}\lambda + c_3 \tilde{ 
          \varepsilon}^{\frac{\gamma_2-1}{\gamma_2}} K^{\sigma} \lambda, \nonumber
\end{eqnarray} 
with constants $c_1 \equiv c_1(n,N,m,\gamma_1,\gamma_2,L/\nu,M)$, 
$c_2 \equiv c_2( n,N,m,\gamma_1,\gamma_2,L/\nu,M,q)$ and $c_3 \equiv \linebreak c_3(n,N,m,\gamma_1,\gamma_2,L/\nu)$.

\subsubsection{Estimate of the maximal function on level sets}\label{chap:Est.Maxfkt}
At this point of the proof we combine the a priori estimate for the solution of the frozen problem with the comparison estimate in order to estimate the super level sets of the maximal function of $|D^m u|^{p_2}$. We use Sobolev-Poincar\'e's inequality to translate the a priori higher differentiability of the solution of the frozen problem into higher integrability and therefore gain an exponent which we denote $n_\delta^*$ (see \eqref{Def.n*delta}). This exponent determines the decay of the super level sets of the maximal function on increasing levels (see \eqref{MF Iteration 2}) and finally provides the desired higher integrability result.

We define the restricted maximal function to the cube $\frac{3}{2}\tilde{Q}$ by
\begin{equation*}
     M^{**} := M^{*}_{\frac{3}{2}\tilde{Q}},
\end{equation*}
whereas 
\begin{equation*}
     M^{*} := M^*_{Q_{4R_0}}
\end{equation*}
denotes the maximal function on $Q_{4R_0}$ (see the statement of Lemma 
\ref{Calderon-Zygmund}).

For $1 < B < B_M$ we now estimate the measure of the set
\begin{equation*}
     \Bigl\{ x \in Q:\ M^{**}(|D^m u|^{p_2})(x) > \tfrac{ABK^{\sigma}\lambda}{2},\
     M^*(|F(\cdot)|^{p(\cdot)}+1)(x) < \tilde{\varepsilon}\lambda \ \Bigr\},
\end{equation*}
where $A$ will be chosen later.
First, by Lemma \ref{techn. 3} we have:
\begin{equation}\label{CZ:splitting}
\begin{aligned}
     |D^m u|^{p_2}
     &\le (\mu^2 + |D^m u|^2)^{p_2/2}\\ 
     &\le c_1 (\mu^2 + |D^m w|^2)^{p_2/2} + c_2 (\mu^2 + |D^m u|^2 + 
           |D^m w|^2)^{\frac{p_2-2}{2}}|D^m u - D^m w|^2\\
     &= \tilde{c}_1{\mathcal G}_1 + \tilde{c}_2{\mathcal G}_2,
\end{aligned}
\end{equation}
where ${\mathcal G}_1 := (\mu^2 + |D^m w|^2)^{p_2/2}$, 
${\mathcal G}_2 := (\mu^2 + |D^m u|^2 
+ |D^m w|^2)^{\frac{p_2-2}{2}}|D^m u - D^m w|^2$ and constants $\tilde{c}_1, 
\tilde{c}_2 \equiv \tilde{c}_1,\tilde{c}_2(n,\gamma_1,\gamma_2)$.
Then there holds
\begin{eqnarray*}
     &&\mskip-20mu \left|\left\{ x \in Q: M^{**} (|D^m u|^{p_2})(x) > 
          \tfrac{ABK^{\sigma}\lambda}{2},\ M^*(|F(\cdot)|^{p(\cdot)}+1)(x) < 
          \tilde{\varepsilon} \lambda \right\}\right|\\
     &&\mskip+30mu \le \left|\left\{ x \in Q: M^{**}({\mathcal G}_1)(x) 
          > \tfrac{ABK^{\sigma}\lambda}{2c_1},\ M^*(|F(\cdot)|^{p(\cdot)}+1)(x) < \tilde{ 
          \varepsilon} \lambda \right\} \right|\\
     &&\mskip+60mu + \left|\left\{ x\in Q: M^{**}({\mathcal G}_2)(x) > 
          \tfrac{ABK^{\sigma}\lambda}{2c_2},\ M^*(|F(\cdot)|^{p(\cdot)}+1)(x) < \tilde{ 
          \varepsilon} \lambda \right\} \right|\\
     &&\mskip+30mu =: I_1 + I_2.
\end{eqnarray*}

\noindent
{\bf Estimate for \boldmath $I_2$\unboldmath }:
Using property (M1) for the maximal function , the inclusion $\frac{3}{2} \tilde{Q} = 
\frac{3}{4} S \subset S$ and the comparison estimate (\ref{Vgl.Absch.}) we obtain
\begin{eqnarray*}
     I_2
     &\le& \frac{c(n)c_2}{ABK^{\sigma} \lambda} \int_S (\mu^2 + |D^m u|^2 +
           |D^m w|^2)^{\frac{p_2-2}{2}}|D^m u - D^m w|^2\, dx\\
     &\le& \frac{c}{ABK^{\sigma}\lambda} \left( c_1 \omega(s) \log \left(\frac{1}{ 
           s}\right) K^{\sigma} \lambda + c_2 \omega(s) \tilde{\sigma}^{-1} K^{\sigma}
           \lambda + c_3 \tilde{\varepsilon}^{\frac{\gamma_2-1}{\gamma_2}} 
           K^{\sigma}\lambda\right) |S|\\
     & = & \frac{\hat{c}_1}{AB} \omega(s) \log \left(\frac{1}{s}\right) |S| + 
           \frac{\hat{c}_2}{AB} \omega(s) \tilde{\sigma}^{-1} |S| + 
           \frac{\hat{c}_3}{AB} \tilde{\varepsilon}^{\frac{\gamma_2-1}{\gamma_2}} |S|,
\end{eqnarray*}
with constants  $\hat{c}_1 \equiv \hat{c}_1(n,N,m,\gamma_1,\gamma_2,L/\nu,M), ,\hat{c}_2 \equiv \hat{c}_2(n,N,m,\gamma_1,\gamma_2,L/\nu,M,q)$ and 
$\hat{c}_3 \equiv \linebreak \hat{c}_3 (n,N,m,\gamma_1,\gamma_2,L/\nu)$.

\noindent
{\bf Estimate for \boldmath $I_1$\unboldmath }:
Since $\frac{n\chi}{n-\chi}$ is increasing in $\chi$, we can find $\delta \equiv
\delta(n,m,\gamma_1,\gamma_2,L/\nu) > 0$, such that
\begin{equation}\label{Ident.delta}
     \frac{n}{n-2}(1+\delta) = \frac{(1+\tilde{\delta})n}{n-2(1+\tilde{\delta})}.
\end{equation}
We set
\begin{equation}\label{Def.n*delta}
     r := \frac{n}{n-2}(1+\delta)p_2 = n^*_{\delta} p_2, \quad \mbox{ with } \
     n^*_{\delta} \equiv \frac{n}{n-2}(1+\delta),
\end{equation}
and distinguish the cases $1 < p_2 < 2$ and $p_2 \ge 2$.

\noindent
{\bf The case \boldmath $p_2 \ge 2$\unboldmath }: 
For $\eta \in \bR^{\cN}$ we estimate with 
\eqref{prop.V} and Lemma \ref{techn. 3} as follows (note that this estimate holds for any 
$p_2 > 1$)
\begin{equation}\label{techn.CZ.1}
\begin{aligned}
     \left(\mu^2 + |D^m w|^2\right)^{p_2/2} 
     &\le c\left( \mu^2 + |\eta|^2\right)^{
          p_2/2} + c\left(\mu^2 + |D^m w|^2 + |\eta|^2\right)^{\frac{p_2-2}{2}}
          |D^m w - \eta|^2\\
     &\le c\left|V_{\mu}(D^m w) - V_{\mu}(\eta)\right|^2 + c\left(\mu^2 + 
          |\eta|^2\right)^{p_2/2},
\end{aligned}
\end{equation}
with the definition of $V_{\mu}(D^m w)$ of \eqref{Def.Vmu,Hmu}.
Therefore by Corollary \ref{Kor.Maxfkt.} we infer that (note also \eqref{Ident.delta})
\begin{eqnarray*}
     I_1 
     &\le& \frac{c(n,r,p_2)c_1^{n^*_{\delta}}}{(ABK^{\sigma} \lambda)^{n^*_{\delta}}} |S|
           \midint_{\frac{3}{4}S} (\mu^2 + |D^m w|^2)^{\frac{p_2}{2} \frac{r}{p_2}}
           \, dx\\
     &\le& \frac{c|S|}{(ABK^{\sigma} \lambda)^{n^*_{\delta}}}
           \midint_{\frac{3}{4}S} \left|V_{\mu}(D^m w) - V_{\mu}(\eta) \right|^{\frac{2n(1+\tilde{\delta})}{n-2(1+\tilde{\delta})}}\, dx
           + \frac{c|S|(\mu^2 + |\eta|^2)^{\frac{p_2 n(1+\delta)}{2(n-2)}}}{(ABK^{\sigma}\lambda)^{n^*_{\delta}}}\\
     & =:& I_{1,1}^{(1)} + I_{1,2}^{(1)},
\end{eqnarray*}
with the obvious labelling of $I_{1,1}^{(1)}$ and $I_{1,2}^{(1)}$. We now choose $\eta$ such that
\begin{equation}\label{choice.eta}
     V_{\mu}(\eta) = \midint_{\frac{3}{4}S} V_{\mu} (D^m w)\, dx.
\end{equation}

\noindent
{\bf Estimate for \boldmath $I_{1,2}^{(1)}$ \unboldmath }: In the case $\mu = 0$ we have the identity
\begin{equation*}
     \left| \left( V_0 (D^m w)\right)_{\frac{3}{4}S} \right| = |V_0(\eta)| = 
     |\eta|^{\frac{p_2-2}{2}}|\eta| = |\eta|^{p_2/2},
\end{equation*}
and therefore together with the energy estimate \eqref{Energ.Absch.} and higher integrability (note that $\lambda \ge 1$, $K \ge 1$) a straight forward estimate
shows
\begin{eqnarray*}
     I_{1,2}^{(1)} 
     &\le& \frac{c|S|}{(ABK^{\sigma}\lambda)^{n^*_{\delta}}} \left( \midint_S |D^m w|^{p_2}\, dx \right)^{n^*_{\delta}}\\
     &\le&\frac{c|S|}{(ABK^{\sigma}\lambda)^{n^*_{\delta}}} \left( \midint_S \left( |D^m u|^{p_2} + 1 \right)
          \, dx \right)^{n^*_{\delta}}
     \le \frac{c|S|}{(AB)^{n^*_{\delta}}}.
\end{eqnarray*}

In the case $\mu \in (0,1]$ we estimate (note that $p_2 \ge 2$):
\begin{equation*}
     |\eta|^{p_2} 
     \le \left|V_{\mu}(\eta)\right|^2 \le \midint_{\frac{3}{4}S} \left|V_{\mu}(D^m w)\right|^2\, dx
     \le \midint_{\frac{3}{4}S} \left(\mu^2 + |D^m w|^2\right)^{p_2/2}\, dx,
\end{equation*}
and therefore with H\"older's inequality, the energy estimate \eqref{Energ.Absch.} and higher
integrability (note that $\mu \le 1$) we easily deduce
\begin{eqnarray*}
     I_{1,2}^{(1)}
     &\le& \frac{c|S|}{(ABK^{\sigma}\lambda)^{n^*_{\delta}}} \left(1+ \midint_{\frac{3}{4}S} \left( \mu^2
           + |D^m w|^2\right)^{p_2/2}\, dx \right)^{n^*_{\delta}}\\
     &\le& \frac{c|S|}{(ABK^{\sigma}\lambda)^{n^*_{\delta}}} \left(\midint_S \left( |D^m u|^{p_2}+1\right)\, dx
           \right)^{n^*_{\delta}}\\
     &\le& \frac{c|S|}{(AB)^{n^*_{\delta}}}.
\end{eqnarray*}

\noindent
{\bf Estimate for \boldmath $I_{1,1}^{(1)}$\unboldmath }:
By the choice of $\eta$ we can apply Sobolev-Poincar\'e's inequality to obtain
\begin{equation*}
     \midint_{\frac{3}{4}S} \left|V_{\mu}(D^m w) - V_{\mu}(\eta) \right|^{\frac{2n(1+\tilde{\delta})}{n-2(1+\tilde{\delta})}}\, dx
     \le c_{SP} \left( s^2 \midint_{\frac{3}{4}S}
     \left|DV_{\mu}(D^m w)\right|^{2(1+\tilde{\delta})}\, dx \right)^{n^*_{\delta}}.
\end{equation*}
The apriori estimate \eqref{Est.w.p_2>2}, taken together with the energy estimate \eqref{Energ.Absch.}
and again higher integrability now provides
\begin{equation*}
     I_{1,1}^{(1)} 
     \le \frac{c|S|}{(ABK^{\sigma} \lambda)^{n^*_{\delta}}} \left( \midint_S 
          \left|H_{\mu}(D^m w)\right|^2\, dx\right)^{n^*_{\delta}}
     \le \frac{c|S|}{(AB)^{n^*_{\delta}}}.
\end{equation*} 
Taking all the estimates together we end up with
\begin{equation*}
     I_1 \le \frac{c|S|}{(AB)^{n^*_{\delta}}},
\end{equation*}
where $c \equiv c(n,N,m,\gamma_1,\gamma_2,L/\nu,M)$.

\noindent
{\bf The case \boldmath $1 < p_2 < 2$\unboldmath }:
For $\mu \in (0,1]$ we first estimate by Corollary \ref{Kor.Maxfkt.} (again note
\eqref{Ident.delta}):
\begin{eqnarray*}
     I_1 
     &\le&\frac{c(n,r,p_2)c_1^{n^*_{\delta}}}{(ABK^{\sigma}\lambda)^{n^*_{\delta}}}|S| \midint_{\frac{3}{4}S}
          H_{\mu}(D^m w)^{\frac{2r}{p_2}}\, dx\\
     &\le&\frac{c|S|}{(ABK^{\sigma}\lambda)^{n^*_{\delta}}} \midint_{\frac{3}{4}S} \left| H_{\mu}(D^m w)
          - \left( H_{\mu}(D^m w)\right)_{\frac{3}{4}S} \right|^{\frac{2n(1+\tilde{\delta})}{n-2(1+\tilde{\delta})}}\, dx
          + \frac{c|S|\left| \left(H_{\mu}(D^m 
          w)\right)_{\frac{3}{4}S} \right|^{n^*_{\delta}}}{(ABK^{\sigma}\lambda)^{n^*_{\delta}}}\\
     & = &I_{1,1}^{(2)} + I_{1,2}^{(2)}.
\end{eqnarray*}
with the definition for $H_{\mu}(D^m w)$ of \eqref{Def.Vmu,Hmu} and the obvious labelling of
$I_{1,1}^{(2)}$ and $I_{1,2}^{(2)}$. 

\noindent
{\bf Estimate for \boldmath $I_{1,1}^{(2)}$\unboldmath }:
Applying Sobolev-Poincar\'e's inequality, the a priori estimate \eqref{Est.w.p_2<2.1} and finally
the energy estimate \eqref{Energ.Absch.} and higher integrability, we obtain
\begin{eqnarray*}
     I_{1,1}^{(2)}
     &\le&\frac{c|S|}{(ABK^{\sigma}\lambda)^{n^*_{\delta}}} \left( s^2 \midint_{\frac{3}{4}S}
          \left|DH_{\mu}(D^m w)\right|^{2(1+\tilde{\delta})}\, dx \right)^{n^*_{\delta}}\\
     &\le&\frac{c|S|}{(ABK^{\sigma}\lambda)^{n^*_{\delta}}} \left( \midint_S \left| H_{\mu}(D^m w)
          \right|^2\, dx \right)^{n^*_{\delta}}\\
     &\le&\frac{c|S|}{(ABK^{\sigma} \lambda)^{n^*_{\delta}}} \left( \midint_S \left( |D^m u|^{p_2}+1
          \right)\, dx \right)^{n^*_{\delta}}
     \le\frac{c|S|}{(AB)^{n^*_{\delta}}}.
\end{eqnarray*}

\noindent
{\bf Estimate for \boldmath $I_{1,2}^{(2)}$\unboldmath }: Here we use H\"older's inequality, the energy
estimate \eqref{Energ.Absch.} and higher integrability to conclude
\begin{eqnarray*}
     I_{1,2}^{(2)}
     &\le&\frac{c|S|}{(ABK^{\sigma} \lambda)^{n^*_{\delta}}} \left( \midint_S \left( |D^m w|^{p_2}+1
          \right)\, dx \right)^{n^*_{\delta}}\\
     &\le&\frac{c|S|}{(ABK^{\sigma} \lambda)^{n^*_{\delta}}} \left( \midint_S|D^m u|^{p_2}\, dx + 1
          \right)^{n^*_{\delta}}\\
     &\le&\frac{c|S|}{(AB)^{n^*_{\delta}}}.
\end{eqnarray*}

In the case $\mu = 0$ we proceed as follows, again using \eqref{techn.CZ.1}:
\begin{equation*}
     |D^m w|^{p_2} \le c \left|V_0(D^m w) - V_0(\eta)\right|^2 + c |\eta|^{p_2}.
\end{equation*}
Therefore we write by Corollary \ref{Kor.Maxfkt.} and again noting \eqref{Ident.delta}:
\begin{eqnarray*}
     I 
     &\le& \frac{c|S|}{(ABK^{\sigma} \lambda)^{n^*_{\delta}}} \midint_{\frac{3}{4}S}|D^m w|^{\frac{np_2(1+\delta)}{
           n-2}}\, dx\\
     &\le& \frac{c|S|}{(ABK^{\sigma} \lambda)^{n^*_{\delta}}} \midint_{\frac{3}{4}S} \left|V_0(D^m w)
           - V_0(\eta)\right|^{\frac{2n(1+\tilde{\delta})}{n-2(1+\tilde{\delta})}}\, dx + \frac{c|S| |\eta|^{\frac{p_2 n(1+\delta)}{n-2}}}{(ABK^{\sigma} 
           \lambda)^{n^*_{\delta}}}\\
     & = & I_{1,1}^{(3)} + I_{1,2}^{(3)}.
\end{eqnarray*}
We choose $\eta$ such that 
\begin{equation*}
     V_0(\eta) = \midint_{\frac{3}{4}S} V_0(D^m w)\, dx.
\end{equation*}

\noindent
{\bf Estimate for \boldmath $I_{1,1}^{(3)}$\unboldmath }: 
Here we use Sobolev-Poincar\'e's inequality, the apriori estimate \eqref{Est.w.p_2<2.2} and finally 
again the energy estimate \eqref{Energ.Absch.} and higher integrability to conclude
\begin{eqnarray*}
     I_{1,1}^{(3)}
     &\le& \frac{c|S|}{(ABK^{\sigma} \lambda)^{n^*_{\delta}}} \left( s^2 \midint_{\frac{3}{4}S}
           \left|DV_0(D^m w)\right|^{2(1+\tilde{\delta})}\, dx \right)^{n^*_{\delta}}\\
     &\le& \frac{c|S|}{(ABK^{\sigma} \lambda)^{n^*_{\delta}}} \left( \midint_S \left|H_0(D^m         
           w)\right|^2\, dx\right)^{n^*_{\delta}}\\
     &\le& \frac{c|S|}{(AB)^{n^*_{\delta}}}.
\end{eqnarray*}

\noindent
{\bf Estimate for \boldmath $I_{1,2}^{(3)}$\unboldmath }: Since we have $|\eta|^{p_2/2} =
| (V_0(D^m w))_{\frac{3}{4}S}|$, we can estimate $I_{1,2}^{(3)}$ in a completely
analogous way as in the case $p_2 \ge 2$ to obtain
\begin{equation*}
     I_{1,2}^{(3)} \le \frac{c|S|}{(AB)^{n^*_{\delta}}}.
\end{equation*}

Thus we have shown in any case (i.e. for any $p_2 > 1$ and for any $\mu \in [0,1]$):
\begin{equation}\label{CZ:gain}
     I_1 \le \frac{\hat{c}_4|S|}{(AB)^{n^*_{\delta}}},
\end{equation}
with $\hat{c}_4 \equiv \hat{c}_4(n,N,m,\gamma_1,\gamma_2,L/\nu)$. Connecting this with the estimate 
valid for $I_2$ we finally arrive at (eventually enlarging the constants by a factor $c(n)$)
\begin{eqnarray*}
     &&|\{ x \in Q: M^{**}(|D^m u|^{p_2})(x) > \tfrac{AB}{2}K^{\sigma}\lambda,\ 
          M^*(|F(\cdot)|^{p(\cdot)}+1)(x) < \tilde{\varepsilon}\lambda\}|\\[0.1cm]
     &&\mskip+100mu \le \Bigl[ \frac{\hat{c}_1}{AB}\omega(s)\log \left(\tfrac{1}{s}
          \right) + \frac{\hat{c}_2}{AB}\omega(s) \tilde{\sigma}^{-1} + 
          \frac{\hat{c}_3}{AB} \tilde{\varepsilon}^{\frac{\gamma_2-1}{\gamma_2}} + \frac{ 
          \hat{c}_4}{(AB)^{n^*_{\delta}}} \Bigr] |Q|.
\end{eqnarray*}

Now we come to the rather involved choice of the parameters. First we determine
$R_1 \equiv R_1(n,N,m,\gamma_1,\gamma_2, L/\nu,\omega(\cdot),\tilde{\sigma}, A,B_M)$ small enough to have
\begin{equation*}
     \frac{\hat{c}_1}{A} \omega(s) \log \left(\tfrac{1}{s}\right) \le \frac{1}{8
     B_M^{n^*_{\delta}-1}} \qquad \mbox{and} \qquad 
     \frac{\hat{c}_2}{A} \omega(s) \tilde{\sigma}^{-1} \le \frac{1}{8B_M^{n^*_{\delta}-1}},
\end{equation*}
for all $s \le 8nR_1$. Then if $R_0 ─\le R_1$ satisfies (\ref{Lokalisierung}) and
(\ref{Einschr.R_0}), we have
\begin{equation*}
R_0 \equiv R_0(n,N,m,\gamma_1,\gamma_2,\nu,L, \| |D^m u(\cdot)|^{p(\cdot)} \|_{L^1}, \| |F(\cdot)|^{p(\cdot)} \|_{L^{\frac{n}{n-\gamma_1}}},\omega(\cdot),B_M).
\end{equation*}
Next we choose $\tilde{\varepsilon} \equiv \tilde{\varepsilon}(n,N,m,\gamma_1,\gamma_2,\nu,L,A,B) \in (0,1)$ such that
\begin{equation}\label{choice.eps}
     \frac{\hat{c}_3}{A} \tilde{\varepsilon}^{\frac{\gamma_2-1}{\gamma_2}} =
     \frac{1}{8 B^{n^*_{\delta}-1}}.
\end{equation}
Next we fix $A$ by
\begin{equation}\label{Def.A}
     A = \max\{ (8 \hat{c}_4)^{n^*_{\delta}}, 5^{n+1}\} \ge 2,
\end{equation}
which yields
\begin{equation*}
     \frac{\hat{c}_4}{(AB)^{n^*_{\delta}}}  \le \frac{1}{8B^{n^*_{\delta}}}.
\end{equation*}
Noting that 
\begin{equation*}
     \frac{1}{8B_M^{n^*_{\delta}-1}} \le \frac{1}{8B^{n^*_{\delta}-1}},
\end{equation*}
for $1 < B \le B_M$ we obtain for any $R \le R_0$ that
\begin{equation*}
     |\{ x \in Q:\ M^{**}(|D^m u|^{p_2})(x) > \tfrac{AB}{2}K^{\sigma}\lambda,\
     M^*(|F(\cdot)|^{p(\cdot)}+1)(x) < \tilde{\varepsilon}\lambda \}| \le 
     \frac{|Q|}{2B^{n^*_{\delta}}}.
\end{equation*}
In particular for every $\vartheta$ satisfying $\vartheta \ge \frac{1}{B^{n^*_{\delta}}}$ there holds
\begin{equation}\label{Maxfkt.delta.1}
     |\{ x \in Q:\ M^{**}(|D^m u|^{p_2})(x) > \tfrac{AB}{2}K^{\sigma}\lambda,\
     M^*(|F(\cdot)|^{p(\cdot)}+1)(x) < \tilde{\varepsilon}\lambda \}| 
     \le \ \tfrac{\vartheta}{2}|Q|.
\end{equation}

We next want to turn this estimate for the maximal function with the fixed
exponent $p_2$ into an estimate for the maximal function of $|D^m u|^{p(\cdot)}$.
Since $p_2 \ge p(x)$ for any $x \in 2 \tilde{Q}$, we see that for any cube
$Q \subset \frac{3}{2}\tilde{Q}$ we have
\begin{equation*}
     \midint_Q |D^m u|^{p(x)}\, dx \le \midint_Q |D^m u|^{p_2}\, dx +1.
\end{equation*}
Hence, for $x \in Q$ there holds
\begin{equation*}
     M^{**}(|D^m u|^{p(\cdot)})(x) \le M^{**}(|D^m u|^{p_2}+1)(x).
\end{equation*}
Since $\lambda,K,\frac{A}{2},B \ge 1$, we have in particular that
$\frac{AB}{2}K^{\sigma}\lambda \ge 1$ and therefore
\begin{equation*}
M^{**}(|D^m u(\cdot)|^{p(\cdot)})(x) > ABK^{\sigma}\lambda 
\end{equation*}
implies
\begin{equation*}
     M^{**}(|D^m u|^{p_2})(x) + \tfrac{AB}{2}K^{\sigma}\lambda 
     \ge M^{**}(|D^m u|^{p_2}+1)(x) > ABK^{\sigma}\lambda
     = 2 \cdot \tfrac{AB}{2}K^{\sigma}\lambda.
\end{equation*}
From \eqref{Maxfkt.delta.1} we therefore obtain
\begin{equation}\label{Maxfkt.delta.2}
     |\{ x \in Q:
     M^{**}(|D^m u(\cdot)|^{p(\cdot)})(x) > ABK^{\sigma}\lambda,
     M^*(|F(\cdot)|^{p(\cdot)}+1)(x) < \tilde{\varepsilon}\lambda\}|\le \tfrac{\vartheta}{2}|Q|.
\end{equation}

The last step in the proof consists in converting \eqref{Maxfkt.delta.2} into an
estimate for the restricted maximal function $M^* = M^*_{Q_{4R_0}}$. This can be
achieved by looking carefully at the cubes involved in the proof. 
Let $\ell$ be the sidelength of the cube $Q$. For
an arbitrary point $x \in Q$ both $x$ itself and the point $x_0$ chosen in 
(\ref{CS 3}) are contained in the cube $\tilde{Q}$ which has sidelength $2\ell$.\\
Now if $C' \subseteq Q_{4R_0}$ is a cube, containing $x$ and having side length 
$\ell'$ larger than $\ell/2$, there holds $C' \cap \tilde{Q} \not= \emptyset$. Thus there exists a cube $C'' \subseteq Q_{4R_0}$, containing $C'$ and $\tilde{Q}$, and whose side length $\ell''$ is bounded by
\begin{equation*}
     \ell'' \le 2\ell+\ell' \le 5\ell'.
\end{equation*}
Therefore, by (\ref{CS 3}) there holds
\begin{equation*}
     \midint_{C'} |D^m u|^{p(x)}\, dx \le \frac{1}{|C'|} \midint_{C''} |D^m u|^{
     p(x)}\, dx \le \frac{|C''|}{|C'|} \lambda \le 5^n \lambda,
\end{equation*}
while in the case $\ell' \le \frac{\ell}{2}$, we have $C' \subset \frac{3}{2} 
\tilde{Q}$ and
\begin{equation*}
     \midint_{C'} |D^m u|^{p(x)}\, dx \le M^{**}(|D^m u|^{p(\cdot)})(x).
\end{equation*}
This implies that
\begin{equation*}
     M^*(|D^m u|^{p(\cdot)})(x) \le \max \left\{ M^{**}(|D^m u|^{p(\cdot)})(x), 5^n 
     \lambda \right\} \quad \mbox{ for all }\ x \in Q.
\end{equation*}
From the choice of $A$, i.e. (\ref{Def.A}) we infer that $\frac{AB}{2}K^{\sigma} \ge \frac{5^{n+1}}{2} > 5^n$.
\begin{equation*}
     \Bigl\{ x \in Q: M^*(|D^m u(\cdot)|^{p(\cdot)})(x) > ABK^{\sigma}\lambda\Bigr\}
     \subseteq \Bigl\{ x \in
     Q: M^{**}(|D^m u(\cdot)|^{p(\cdot)})(x) > ABK^{\sigma}\lambda\Bigr\},
\end{equation*}
and therefore
\begin{equation*}
     |\{ x \in Q:
     M^{*}(|D^m u(\cdot)|^{p(\cdot)})(x) > ABK^{\sigma}\lambda,
     M^*(|F(\cdot)|^{p(\cdot)}+1)(x) < \tilde{\varepsilon}\lambda \}|
     \le \tfrac{\vartheta}{2}|Q|.
\end{equation*}
This contradicts (\ref{Calderon-Zygmund 1}) and completes the proof of Lemma 
\ref{Calderon-Zygmund}.
\end{proof}

\subsection{Proof of the main theorem}\label{ch:CZ.main}

We now apply Lemma \ref{Calderon-Zygmund} in order to obtain the result of the main theorem. Since the procedure follows the one of \cite[pp 141-146]{Acerbi:2004}, we only sketch the main steps here, nevertheless explicitely pointing out the special choice of the constants and parameters.

First we define 
\begin{eqnarray*}
     &&\mu_1(t) := \left|\left\{ x \in Q_{R_0}: M^*(|D^m u(\cdot)|^{p(\cdot)})(x) > 
          t\right\}\right|,\\[0.2cm]
     &&\mu_2(t) := \left|\left\{ x \in Q_{R_0}: M^*_{1+\sigma} 
          (|F(\cdot)|^{p(\cdot)}+1)(x) > t\right\} \right|,
\end{eqnarray*}
with $M^* \equiv M^*_{Q_{4R_0}}$. For $q \in (1,\frac{n}{n-2})$ we set
\begin{equation}\label{Def. B_M}
     B_M := \left( 2(AK_M^{\sigma_M})^q\right)^{\frac{n-2}{n(1+\delta)-q(n-2)}}
     = B_M(n,q,\delta,K_M,\sigma_M).
\end{equation}
By the restriction imposed on the range of $q$ there holds $\frac{n-2}{n(1+\delta)-q(n-2)} >
\frac{n-2}{n\delta + 2} > 0$ and therefore, since $A \ge 2, K_M > 1, \sigma_M > 0$, we have $B_M > 1$.
With this choice of $B_M$ we set for $1 < K < K_M$ and $0 < \sigma < \sigma_M$:
\begin{equation}\label{Def. B}
     B := \left( 2(AK^{\sigma})^q\right)^{\frac{n-2}{n(1+\delta)-q(n-2)}}.
\end{equation}
Then
\begin{equation*}
     1 < B < B_M,
\end{equation*}
and moreover
\begin{equation}\label{Bez. B,AK}
B^{-n^*_{\delta} + q} = \frac{1}{2(AK^{\sigma})^q}.
\end{equation}
Now, we let
\begin{equation}\label{Def. lambda_0}
     \lambda_0 := \frac{5^{n+2}c_W}{\vartheta} \midint_{Q_{4R_0}} |D^m u|^{p(x)}\, dx +1,
\end{equation}
where we have chosen
\begin{equation}\label{Wahl delta}
\vartheta := B^{-n^*_{\delta}}.
\end{equation}
Here $c_W \equiv c_W(n)$ denotes the constant from (\ref{Absch.Mfkt.1}).
By (\ref{Absch.Mfkt.1}) and the definition of $\lambda_0$ we obtain
\begin{equation}\label{Absch.mu1lambda0}
     \mu_1(\lambda_0) 
     \le \frac{c_W}{\lambda_0} \midint_{Q_{4R_0}} |D^m u|^{p(x)}\, dx
           \cdot |Q_{4R_0}|
     \le \frac{4^n |Q_{R_0}|\vartheta}{5^{n+2}} \le \tfrac{\vartheta}{2} |Q_{R_0}|.
\end{equation}
Let $A$ be the constant from Lemma \ref{Calderon-Zygmund}. 
Since $A,B,K \ge 1$ we have 
$ABK^{\sigma} \ge 1$ and therefore
\begin{equation}\label{Absch mu_1}
     \mu_1\left((ABK^{\sigma})^h \lambda_0\right) \le \tfrac{\vartheta}{2} |Q_{R_0}|
     \qquad \mbox{for}\quad h \in \bN \cup \{0\}.
\end{equation}
We now set 
\begin{equation}\label{Def.tildeA}
\tilde{A} := ABK^{\sigma} \ge 2.
\end{equation}
It can easily be checked that with $h \in \bN \cup \{0\}$ the assumptions of Lemma \ref{Calderon-Zygmund} are satisfied for the sets 
$$X := \{ x \in Q_{R_0}: M^*(|D^m u(\cdot)|^{p(\cdot)})(x) > 
          \tilde{A}^{h+1}\lambda_0,
          M^*_{1+\sigma} (|F(\cdot)|^{p(\cdot)}+1)(x) < \tilde{\varepsilon}(AB)^h 
          \lambda_0\}
$$ 
and 
$$Y := \{ x \in Q_{R_0}:\ M^*(|D^m u(\cdot)|^{p(\cdot)})(x) > \tilde{A}^h\lambda_0\},
$$
which provides the estimate
\begin{equation}\label{MF Iteration 1}
     \mu_1( \tilde{A}^{h+1}\lambda_0) \le B^{-n^*_{\delta}}
     \mu_1( \tilde{A}^h \lambda_0) + \mu_2(\tilde{\varepsilon}
     \tilde{A}^h\lambda_0)
\end{equation}
Iterating this inequality and exploiting the definition of $A$ and the specific choice of $B$ we obtain for arbitrary $J \in \bN$
\begin{equation}\label{MF Iteration 2}
     \mu_1(\tilde{A}^{h+1} \lambda_0) \le \mu_1(\lambda_0)
     + \tilde{A}^q \sum\limits_{i=0}^J \mu_2(\tilde{\varepsilon} 
     \tilde{A}^i \lambda_0) \tilde{A}^i \sum\limits_{k=0}^{J-i}\tilde{A}^{qk} 
     B^{-n^*_{\delta}k}.
\end{equation}

\noindent
Using the choice of $B$ and the definition of $\tilde{A}$ we infer that 
     $\tilde{A}^{qk} B^{-n^*_{\delta}k} = ( (AK^{\sigma})^q B^{q-n^*_{\delta}})^k = 2^{-k}$.
Hence, the last sum on the right hand side above can be uniformly estimated from above by 
\begin{equation*}
     \sum\limits_{k=0}^{J-i} \tilde{A}^{qk} B^{n^*_{\delta}k} \le 2.
\end{equation*}
Since the estimate holds for any $J \in \bN$, we obtain, passing to the limit $J \to \infty$:
\begin{equation}\label{MF Absch. 1}
     \sum\limits_{k=1}^{\infty} \tilde{A}^{qk} \mu_1( \tilde{A}^k
     \lambda_0)\le \mu_1(\lambda_0) + 2\tilde{A}^q\sum\limits_{k=0}^{\infty} 
     \tilde{A}^{kq} \mu_2( \tilde{A}^k \lambda_0 \tilde{\varepsilon}).
\end{equation}

This estimate can be transformed into an estimate for the maximal function. Applying the elementary identity
\begin{equation*}
     \int_Q g^q\, dx = \int_0^{\infty} q\lambda^{q-1}\left| \{x \in Q:\ g(x) > \lambda \}\right|\, 
     d\lambda,
\end{equation*}
which holds for $g \in L^q(Q), g \ge 0$, $q \ge 1$, to the maximal function of $|D^m u|^{p(\cdot)}$, decomposing the inteval $[0,\infty)$ into intervals $[0,\lambda_0]$ and $[\tilde{A}^n\lambda_0,\tilde{A}^{n+1}\lambda_0]$ and exploiting \eqref{MF Absch. 1} in combination with the monotonicity of $\mu_1(t)$ and $\mu_2(t)$ and finally using the $L^p$ estimate for the maximal function, we end up with the estimate
\begin{equation}\label{CZ:est.right.3}
     \int_{Q_{R_0}} |D^m u|^{p(x)q}\, dx
     \le |Q_{R_0}|\lambda_0^q + 2(\tilde{A}\lambda_0)^q \mu_1(\lambda_0)
          + c(n) \frac{q^2}{q-1} \cdot \frac{\tilde{A}^{2q}}{\tilde{
          \varepsilon}^q} \int_{Q_{4R_0}} \mskip-10mu |F|^{p(x)q}+1\, dx.
\end{equation}
By (\ref{Absch.mu1lambda0}), by the choice of $\vartheta$ in (\ref{Wahl delta}) and
\eqref{Bez. B,AK} we obtain
\begin{equation*}
     2(\tilde{A}\lambda_0)^q \mu_1(\lambda_0)
     \le (\tilde{A}\lambda_0)^q B^{-n^*_{\delta}} |Q_{R_0}|
     \overset{(\ref{Bez. B,AK})}{=}  (AK^{\sigma})^q \lambda_0^q 
          \frac{1}{2(AK^{\sigma})^q} |Q_{R_0}|
     = \tfrac{1}{2} \lambda_0^q |Q_{R_0}|.
\end{equation*}
Furthermore, recalling the definition of $\lambda_0$ in \eqref{Def. lambda_0} and the choice of $\vartheta$, and taking into account (\ref{Def. B}), as well as the dependencies of the constant $A$ (recall also the definition of $n^*_{\delta}$ in
\eqref{Def.n*delta}, we see that
\begin{equation*}
     \lambda_0 
     = c(n) B^{n^*_{\delta}} \midint_{Q_{4R_0}} \left(|D^m u|^{p(x)}+1
           \right)\, dx
     = cK^{\frac{\sigma n(1+\delta)q}{n(1+\delta)-q(n-2)}}
           \midint_{Q_{4R_0}} \left(|D^m u|^{p(x)}+1\right)\, dx,
\end{equation*}
with $c \equiv c(n,\gamma_1,\gamma_2,\nu,L,M)$.
To treat the last term appearing on the right hand side of \eqref{CZ:est.right.3}
we need to control $\tilde{A}^{2q}/\tilde{\varepsilon}^q$. Recalling the definitions
of $\tilde{\varepsilon}, \tilde{A}$ and $B$ from \eqref{choice.eps}, \eqref{Def.tildeA} and \eqref{Def. B} we find that
\begin{equation*}
     \frac{\tilde{A}^2}{\tilde{\varepsilon}}
     = (ABK^{\sigma})^2 \left(\frac{A}{8B^{\frac{n(1+\delta)}{n-2}-1}} \cdot \frac{1}{
           \hat{c}_3} \right)^{-\frac{\gamma_1}{\gamma_1-1}}
     = c K^{\frac{2q\sigma(n-2)}{n(1+\delta)-q(n-2)} + \frac{2q\gamma_1 \sigma}{(n(1+\delta)-q(n-2))( 
           \gamma_1-1)} + 2\sigma},
\end{equation*}
where $c \equiv c(n,\gamma_1,\gamma_2,L/\nu,M,q)$. 
Thus passing to the averages leeds to
\begin{eqnarray*}
     \left(\midint_{Q_{R_0}} |D^m u|^{p(x)q}\, dx \right)^{1/q}
     &\le&cK^{\frac{\sigma n(1+\delta)q}{n(1+\delta)-q(n-2)}} \midint_{Q_{4R_0}} |D^m 
           u|^{p(x)}+1\, dx\\
     &&+ cK^{\frac{2(n-2)q\sigma}{n(1+\delta)-q(n-2)} + \frac{2q\gamma_1 
           \sigma}{(n(1+\delta)-q(n-2))(\gamma_1-1)} + 2\sigma} \left(\midint_{Q_{4R_0}} 
           \mskip-10mu |F|^{p(x)q}+1\, dx\right)^{1/q}.
\end{eqnarray*}

For given $\varepsilon > 0$, we
now want to reach the following smallness conditions to be fulfilled:
\begin{equation}\label{Bed.sigma.1}
     \frac{\sigma n(1+\delta) q}{n(1+\delta)-q(n-2)} \le \varepsilon
\end{equation}
and
\begin{equation}\label{Bed.sigma.2}
     \frac{2(n-2)\sigma q}{n(1+\delta)-q(n-2)} \le \frac{\varepsilon}{3},\quad
     \frac{2\gamma_1\sigma q}{(n(1+\delta)-q(n-2))(\gamma_1-1)} \le \frac{\varepsilon}{3}, \quad
     2\sigma \le \frac{\varepsilon}{3}.
\end{equation}
These conditions hold for example, if
\begin{equation}\label{Wahl sigma}
     \sigma \le \frac{\varepsilon}{9}\min\left\{ 1, \frac{(n(1+\delta)-q(n-2))(\gamma_1-1)
     }{q\gamma_1}, \frac{n(1+\delta)-q(n-2)}{n(1+\delta)q}\right\}.
\end{equation}
Since we have (choosing $\delta$ small enough)
     $\frac{(n(1+\delta)-q(n-2))(\gamma_1-1)}{q\gamma_1} < \frac{n(1+\delta)-q(n-2)}{q} < \frac{n(1+\delta)}{q}-n+2
     < 3$, and $\frac{n(1+\delta)-q(n-2)}{nq} < \frac{2}{n} + \delta < 3$, \eqref{Wahl sigma} implies
     $\sigma < \frac{\varepsilon}{3}$.
To reach this, we set
\begin{equation*}
     \bar{\sigma} := \frac{\varepsilon}{9\sigma_M} \min \left\{ 1, 
     \frac{(n(1+\delta)-q(n-2))(\gamma_1-1)}{q\gamma_1}, \frac{n(1+\delta)-q(n-2)}{n(1+\delta)q} \right\}.
\end{equation*}
and finally
\begin{equation}\label{Wahl tildesigma}
     \tilde{\sigma} = \min \left\{ \bar{\sigma},\gamma_1-1,\tfrac{1}{2} \right\}.
\end{equation}
Thus $\tilde{\sigma} \equiv \tilde{\sigma}(n,q,\gamma_1,\varepsilon,c_g,\delta)$ is fixed.
With this choice, there holds (note that $\sigma_0 < \sigma_M$ and the
estimates above):
\begin{equation*}
     \sigma = \tilde{\sigma} \sigma_0 = \frac{\varepsilon \sigma_0}{9\sigma_M}
     \min\left\{1, \frac{(n(1+\delta)-q(n-2))(\gamma_1-1)}{q\gamma_1}, \frac{n(1+\delta)-q(n-2)}{n(1+\delta)q} 
     \right\} < \frac{\varepsilon}{3} < \varepsilon.
\end{equation*}
Then $\sigma \le \varepsilon$ and $\varepsilon < q-1$ implies (note also that 
$|Q_{4R_0}| \le 1$ since $8nR_0 \le 1$):
\begin{equation*}
     K = \int_{Q_{4R_0}} \left( |D^m u|^{p(x)} + |F|^{p(x)(1+\sigma)}\right)\, dx + 1
     \le \int_{Q_{4R_0}}\left(|D^m u|^{p(x)} + |F|^{p(x)(1+\varepsilon)}\right)\, dx 
     + 2.
\end{equation*}

\noindent
{\bf Remark on the dependencies of the constants:}
By the choice of $\tilde{\sigma} = \tilde{\sigma}(\varepsilon)$,
also $R_1 \equiv R_1(\tilde{\sigma}) = R_1(\varepsilon)$ is fixed via Lemma \ref{Calderon-Zygmund}, and then also $R_0 \equiv R_0(\varepsilon)$ via
Lemma \ref{Calderon-Zygmund} and \eqref{Einschr.R_0}. 

Therefore for any cube $Q_R$ with $R \le R_0$, $Q_{4R} \Subset \Omega$ there holds
\begin{equation*}
     \left(\midint_{Q_{R}}\mskip-12mu |D^m u|^{p(x)q}\, dx \right)^{1/q}\mskip-8mu 
     \le cK^{\varepsilon} \midint_{Q_{4R}}\mskip-12mu |D^m u|^{p(x)}+1\, dx + 
     cK^{\varepsilon} \left(\midint_{ Q_{4R}} \mskip-12mu |F|^{p(x)q}+1\, 
     dx\right)^{1/q},
\end{equation*}
in which the constant depends on $n,N,m,\gamma_1,\gamma_2,L/\nu$ and $q$, and with
\begin{equation*}
     K = \int_{Q_{4R}}\left(|D^m u|^{p(x)} + |F|^{p(x)(1+\varepsilon)}\right)\, dx 
     + 1.
\end{equation*}
Therefore the statement $|D^m u|^{p(\cdot)} \in L^q_{loc}(\Omega)$ follows by a
covering argument. \hfill $\boxed{QED}$


\bibliographystyle{plain}
\label{bibliography}
\makeatletter
\addcontentsline{toc}{chapter}{\bibname}
\makeatother
\bibliography{references}
\nocite{Habermann:Diss, Kristensen:2006, Kronz:Habil, Campanato:1983, Coscia:1999,Manfredi:PHD,Marcellini:1991,Rajagopal:2001,Rao:1991}

\end{document}